\documentclass[oneside,11pt]{article}
\usepackage{stackrel}
\usepackage{graphicx, subfigure}
\usepackage[top=1in,bottom=1in,left=1in,right=1in]{geometry}
\usepackage[sort&compress]{natbib} 
\usepackage{amsfonts, amsmath, amssymb, amsthm, constants, bbm}
\usepackage{MnSymbol}
\usepackage{mathtools}
\usepackage{enumitem}
\usepackage{caption}
\usepackage{algpseudocode}
\usepackage{algorithm2e}
\RestyleAlgo{ruled}
\renewcommand{\leq}{\leqslant}
\renewcommand{\geq}{\geqslant}
\renewcommand{\propto}{\asymp}

\usepackage{cellspace}
\usepackage{makecell}
\setcellgapes{4pt}

\usepackage{color}
\definecolor{darkred}{RGB}{100,0,0}
\definecolor{darkgreen}{RGB}{0,100,0}
\definecolor{darkblue}{RGB}{0,0,150}
\definecolor{citecol}{RGB}{30,80,150}
\definecolor{tabcol}{RGB}{200,230,255}

\usepackage{colortbl}

\usepackage{hyperref}
\usepackage{cleveref}[2012/02/15]
\crefformat{footnote}{#2\footnotemark[#1]#3}
\usepackage[bottom]{footmisc}
\hypersetup{colorlinks=true, linkcolor=darkred, citecolor=citecol, urlcolor=darkblue}
\usepackage{url}

\theoremstyle{theorem}
\newtheorem{thm}{Theorem}
\newtheorem{prp}[thm]{Proposition}
\newtheorem{lem}[thm]{Lemma}
\newtheorem{cor}[thm]{Corollary}

\newtheorem{ass}[thm]{Assumption}

\theoremstyle{remark}
\newtheorem{rem}[thm]{Remark}

\def\beq{\begin{equation}} 
\def\eeq{\end{equation}}
\def\beqn{\begin{eqnarray*}}
\def\eeqn{\end{eqnarray*}}
\def\Bal{\begin{align}}
\def\Eal{\end{align}}
\def\Bitem{\begin{itemize}\setlength{\itemsep}{.2in}}
\def\bitem{\begin{itemize}\setlength{\itemsep}{.05in}}
\def\eitem{\end{itemize}}
\def\blatin{\begin{enumerate}\setlength{\itemsep}{.05in}\renewcommand{\labelenumi}{\roman{enumi}.}}
\def\elatin{\end{enumerate}}
\def\Benum{\begin{enumerate}\setlength{\itemsep}{.2in}}
\def\benum{\begin{enumerate}\setlength{\itemsep}{.05in}}
\def\eenum{\end{enumerate}}
\def\bmult{\begin{multline*}}
\def\emult{\end{multline*}}
\def\bcenter{\begin{center}}
\def\ecenter{\end{center}}
\def\bframe{\begin{frame}}
\def\eframe{\end{frame}}

\newcommand{\thmref}[1]{Theorem~\ref{thm:#1}}
\newcommand{\prpref}[1]{Proposition~\ref{prp:#1}}
\newcommand{\corref}[1]{Corollary~\ref{cor:#1}}
\newcommand{\lemref}[1]{Lemma~\ref{lem:#1}}
\newcommand{\secref}[1]{Section~\ref{sec:#1}}
\newcommand{\figref}[1]{Figure~\ref{fig:#1}}

\newcommand{\assref}[1]{Assumption~\ref{ass:#1}}
\newcommand{\appref}[1]{Appendix~\ref{app:#1}}

\DeclareMathOperator*{\argmax}{arg\, max}

\def\cA{\mathcal{A}}
\def\cB{\mathcal{B}}
\def\cC{\mathcal{C}}
\def\cD{\mathcal{D}}
\def\cE{\mathcal{E}}
\def\cF{\mathcal{F}}
\def\cG{\mathcal{G}}

\def\cI{\mathcal{I}}

\def\cK{\mathcal{K}}
\def\cL{\mathcal{L}}

\def\cN{\mathcal{N}}

\def\cP{\mathcal{P}}

\def\cR{\mathcal{R}}
\def\cS{\mathcal{S}}
\def\cT{\mathcal{T}}
\def\cU{\mathcal{U}}

\def\cX{\mathcal{X}}

\def\bbE{\mathbb{E}}

\def\bbP{\mathbb{P}}

\def\bbR{\mathbb{R}}

\def\ind{\mathbbm{1}}

\newcommand{\E}{\operatorname{\mathbb{E}}}
\renewcommand{\P}{\operatorname{\mathbb{P}}}

\def\Unif{\text{Unif}}

\let\lac\{
\let\rac\}
\renewcommand{\{}{\left\lac}
\renewcommand{\}}{\right\rac}
\newcommand{\inner}[2]{\langle #1, #2 \rangle}

\newcommand{\floor}[1]{\lfloor #1 \rfloor}

\def\1{\mathbbm{1}}

\def\({\left(}
\def\){\right)}

\DeclareMathOperator{\Card}{Card}

\DeclareMathOperator{\Id}{Id}

\newcommand{\ve}{\varepsilon}
\newcommand{\vp}{\varphi}
\newcommand*\diff{\mathop{}\!\mathrm{d}}

\newcommand\wh{\widehat}
\newcommand\dt{\mathrm{d}}

\newcommand\Conv{\mathrm{Conv}}

\definecolor{purple}{rgb}{0.4,.1,.9}

\newcommand{\newclem}[1]{{\color{black}{#1}}}

\newcommand{\newnico}[1]{{\color{black} #1}}

\newcommand{\algoname}{\hyperref[algo:generic]{\texttt{PINES}}\xspace}

\title{Seriation of T\oe plitz and latent position matrices: optimal rates and computational trade-offs}
\author{
	Cl\'ement Berenfeld\footnote{INRIA, PreMeDICaL Team, Univ. Montpellier, France}
	\and
	Alexandra Carpentier\footnote{Institut für Mathematik -- Universität Potsdam, Potsdam, Germany}
	\and
	Nicolas Verzelen\footnote{INRAE, Institut Agro, MISTEA, Univ. Montpellier, France}
}

\begin{document}

\thispagestyle{empty}

\maketitle


\begin{abstract}
In this paper, we consider the problem of seriation of a permuted structured matrix 
based on noisy observations. The entries of the matrix 
relate to an expected quantification of interaction between two objects: the higher the value, the closer the objects. A popular structured class for modeling such matrices 
is the permuted
Robinson class, namely the set of matrices 
whose coefficients are decreasing away from its diagonal, 
up to a permutation of its rows and columns. 
We consider in this paper two submodels of Robinson matrices: the T\oe plitz model, and the latent position model. 
We provide a computational lower bound based on the low-degree paradigm, which hints that there is a statistical-computational gap for
seriation when measuring the error based on the Frobenius norm. We also provide a simple and polynomial-time algorithm that achieves this lower bound. Along the way, we also characterize the information-theory optimal risk thereby giving evidence for the extent of the computation/information gap for this problem. 
\end{abstract}

\section{Introduction}

\subsection{Context and motivation}

The seriation problem is that of ordering $n$ objects from pairwise measurements. Since its introduction in archeology for the chronological dating of  graves~\cite{Robinson51}, it has arisen in various modern data science problems, such as envelope reduction for sparse matrices \cite{barnard1995spectral}, alignment of reads in {\it de novo} sequencing \cite{garriga2011banded,bioinfo17},  time synchronization in distributed networks \cite{Clock-Synchro04,Clock-Synchro06}, or interval graph identification \cite{fulkerson1965incidence}.

In this paper, we consider a setting where we have noisy observations $Y= X+ E$ of the pairwise symmetric interaction matrix $X\in [0,A]^{n\times n}$ for some $A>0$. In general, the noisy seriation problem amounts to recovering a permutation matrix $\Pi^*$ such that that the {\it permuted signal matrix} $\Pi^*X \Pi^{*\top}${\it is a Robinson matrix}, that is the entries of the  $\Pi^*X \Pi^{*\top}$ are non-increasing when one moves away from the diagonal --- see Section~\ref{sec:models} for precise definitions. Such a matrix $X$ is said to be \emph{pre-Robinson}. This property models the fact that objects that are close to each other (with respect to the ordering) tend to have high interactions whereas objects that are farther away tend to have low interactions.  
In this manuscript,  we focus on a specific instance of the seriation problem where  $\Pi^*X \Pi^{*\top}$ also satisfies some stationarity property; namely, we assume that $\Pi^*X \Pi^{*\top}$ is either a {\it T\oe plitz matrix} or has been sampled according to a {\it latent position model} --- see Section~\ref{sec:models}  and the introduction of~\cite{cai2023matrix} for practical motivations e.g. in genomics.

\subsection{Related works}

The seriation problem has attracted a lot of attention both in the computer science literature where the noise matrix $E$, if not null, is sometimes considered as arbitrary and in the statistical literature where $E$ is assumed to have been sampled from some distribution. One interesting feature of the seriation problem is that is exhibits both algorithmic and statistical challenges. 

\paragraph{Robust seriation  with adversarial errors.} There exists an extensive literature in theoretical computer science that aims at recovering the permutation for general Robinson matrices. In the noiseless case ($E=0$), \cite{atkins1998spectral} have established that a simple spectral algorithm is able to recover the permutation $\Pi^*$. See also~\cite{fogel2013convex} for other convex relaxations.  More recently,~\cite{carmona2024modules} have introduced  procedures being able to recover the permutation with an optimal $O(n^2)$ computational complexity. When the perturbation $E$ is deterministic and possibly arbitrary, Chepoi et al. \cite{chepoi2009seriation} have shown that it is NP-hard to recover a pre-Robinson matrix $M'$ such that $\|M'-Y\|_{\infty}\leq \|E\|_{\infty}$, where $\|A\|_{\infty}:= \max_{i,j}|A_{i,j}|$. Conversely, Chepoi and Seston~\cite{chepoi2011seriation} introduced a polynomial-time 16-approximation to that problem, that is they are able to find a pre-Robinson matrix $M'$ such that $\|M'-Y\|_{\infty}\leq 16\|E\|_{\infty}$. For other distances such as the Frobenius distance $\|.\|_F$, the problem is also known to be NP-hard~\cite{barthelemy2001np} and we are not aware of any approximation scheme. While these results are interesting in their own, they do not allow to directly characterize the seriation problem with stochastic noise as NP-hardness results are established for adversarial errors.

\paragraph{Noisy seriation.} Noisy seriation problems have recently gained interest~\cite{flammarion2019optimal,janssen2022reconstruction,giraud2023localization,natik2021consistency,cai2023matrix,issartel2024minimax}. In particular, Flammarion et al.~\cite{flammarion2019optimal} have considered  a related rectangular model where, up to a permutation of the rows, each column is unimodal. The authors have characterized the minimax risk for estimating this matrix, but their procedure, based on a least-square type criterion unfortunately suffers from a high computational cost. In addition, they do not provide any polynomial-time algorithm for this problem. This is in sharp contrast with other shape constraints, such as isotony of columns, which arise in ranking problems, where some polynomial-time procedures provably achieve the minimax risk~\cite{pilliat2024optimal}.
Cai and Ma~\cite{cai2023matrix} consider, as in this manuscript, the noisy permuted T\oe plitz-Robinson matrix. However, they focus on the problem of exactly  recovering the unknown permutation $\Pi^*$ from the observations matrix $Y$.
In~\cite{cai2023matrix}, they characterize the minimum conditions on the vector $\theta$ so that exact reconstruction of $\Pi^*$ is information-theoretically possible. They also establish that the spectral seriation algorithm~\cite{atkins1998spectral} recovers $\Pi^*$ under a much stronger condition, which may lead to conjecture the existence of a computation-statistical gap for this problem, but they did not provide formal evidence for this gap. Moreover, the objective of exactly reconstructing the permutation $\Pi^*$ is much stronger than ours. In particular, when the vector $\theta$ lies in $[0,1]$, exact seriation is possible only when the subgaussian norms of the independent entries of $E$ is bounded by $c/\sqrt{n\log(n)}$, which corresponds to an extremely low-noise situation.

\paragraph{Noisy seriation under additional conditions.} As alluded above, there does not exist polynomial-time procedure with strong theoretical guarantees for general pre-Robinson matrices or for T\oe plitz Robinson matrices. However, there exists a stream of literature (e.g.~\cite{issartel2024minimax,natik2021consistency}) in statistics and machine learning, where the authors put additional conditions on the matrix $X$ and, under this assumption, introduce and analyze polynomial-time  seriation procedures whose error turns out to be minimax optimal. For instance,~\cite{issartel2024minimax} assume that  the entries of $\Pi^* X \Pi^{*\top}$ are bi-Lipschitz. \newnico{This structural assumption is crucial in their work. As revealed by our results, when this structural assumption is removed and we allow for arbitrary T\oe plitz Robinson matrices, a computation-statistical gap occurs.}
Similarly, analysis of the spectral algorithm in~\cite{recanati2018reconstructing,natik2021consistency} make strong assumptions on the matrix $\Pi^* X \Pi^{*\top}$. 
There is also a line of research for seriation problems on graphs. The latter problem is sometimes referred as line embedding on graphons~\cite{janssen2022reconstruction}. As above, available polynomial-time procedures rely on strong additional assumptions on the graphons, although those are not directly expressed in terms of regularity.  

In summary, \newnico{if we make strong structural assumptions on the matrix (as e.g. in~\cite{issartel2024minimax}), the optimal seriation seriation risk is well understood and is achievable in polynomial time}.  \newnico{For general T\oe plitz-Robinson matrices},  one may conjecture from the literature that statistical/computational gaps occur for the seriation problem. 
However, the optimal polynomial-time risks (and even the information-theoretical optimal risks) remain largely unknown, whether for general Robinson matrices or for Toeplitz matrices, despite a significant literature on the topic. 

\subsection{Contributions}

In this work, we measure the quality of a seriation estimator $\hat{\Pi}$ through the so-called $\ell_2$ seriation error defined by 
 $$\ell_2(\hat \Pi):=\inf_{R \in \cR_n}\|\hat \Pi X \hat \Pi^\top -R\|_F ,$$
where $\cR_n$ is the set of Robinson matrices of size $n$ and $\|.\|_F$ is the Frobenius norm. This loss quantifies, in Frobenius distance, to what extent the {\it ordered signal matrix} $\hat \Pi X \hat \Pi^\top$ is close to a Robinson matrix. We further explain the rationale behind this loss in Section~\ref{sec:models}. 
Our contribution is threefold.
\begin{enumerate}
    \item We propose a simple and polynomial-time algorithm \algoname whose $\ell_2$ risk is uniformly bounded in expectation by $n^{3/4}$ up to poly-log terms both for the T\oe plitz and the latent position models.
    \item We provide a matching computational lower bound, in the low-degree polynomial framework~\cite{KuniskyWeinBandeira} suggesting that the rate $n^{3/4}$ cannot be improved by polynomial-time procedures for latent position models. We also provide a similar result for a variant of the T\oe plitz model. \newnico{
 From a broad perspective,  the proof of the computational lower bound focuses on Robinson matrices such that  $\Pi^*X \Pi^{*T}$ is close to a banded matrix with a band of the order of $\sqrt{n}$. Note that such matrices $\Pi^*X \Pi^{*T}$ highly differ from 
 Bi-Lipschitz ones that are considered in~\cite{issartel2024minimax} and for which no computation-statistical gap occurs. More precisely, in our permuted near-banded matrix $X$, on each row, less than $n^{1/2}$ entries bring information on the relative ordering. Our constructions shares some similarities with recent work of Luo and Gao~\cite{luo2023computational} which states computational lower bound for stochastic block models with many groups --- see the discussion section for further details.
    }

    \item Finally, we establish that that the information-theoretic minimax $\ell_2$ risk for this problem is of the order of $n^{1/2}$, thereby providing evidence for the  extent of the statistical-computational gap.  
\end{enumerate}
We also informally extend  \algoname to deal with missing data and consider variations of the procedures to handle seriation for entry-wise errors.  

The idea behind our procedure is to estimate by $\hat{d}$ a suitable distance $\dt$ between two rows of $X$. This distance must have the properties that the \emph{true} neighbors of a row $i$  according to $\Pi^*$   are close to $i$ according to the distance $d$. For any row $i$, we then estimate of set of neighboring rows that are to be  removed from the matrix. Relying on this submatrix, we construct  a geometric graph by connecting two nodes $j$  and $j'$ if and only if their estimated distance is small enough. A key result for this procedure is that, under the above assumptions and if it is properly calibrated, then, for any $i$, this graph has a at most two connected components. These two connected components are made of of $j$'s that are all on left of $i$, and of $j$'s that are right of $i$ with respect to $\Pi^*$. Combining the information of these components for each $i$, we recover a seriation of $X$ that turns out to have the desired properties. While the idea of computing a proxy for the distance is not new --- see e.g.~\cite{issartel2024minimax}--- all previous procedures and analyses rely on specific assumptions on the matrix $\Pi^*X \Pi^{*\top}$. 
Up to our knowledge, this is the first polynomial-time procedure whose error is uniformly bounded on the whole collection of T\oe plitz Robinson matrices or on latent space matrices. We further compare our results to the literature in the Discussion section. 

\subsection{Organization of the paper}

We introduce both the T\oe plitz and the latent position models as well as the corresponding loss functions in Section~\ref{sec:models}. In Section~\ref{sec:l2seriation}, we build our polynomial-time seriation estimator \algoname and we provide uniform risk bounds for both these models. A computational lower bound is provided in Section~\ref{sec:lower:bounds:low_degee} thereby showing the optimality of \algoname. In Section~\ref{sec:noncomp}, we characterize the minimax risks and thereby establish that, if we allow for exponential-time procedures, the risk is significantly smaller. Finally, we further extend our methodology and discuss the literature in Section~\ref{sec:discussion}. As a byproduct of our generic procedure, we also consider other seriation problems in that section. All the proofs are postponed to the appendix.

\paragraph{Notation}

We let $[n]$ be the set of natural numbers between $1$ and $n$. A vector $x \in \bbR^n$ is said to be unimodal with respect to an index $i \in [n]$ if
$$
x_1 \leq \dots \leq x_{i-1} \leq x_i \geq x_{i+1} \geq \dots \geq x_n. 
$$
Likewise, we say the a function $\vp$ defined on a subset $\cU$ of $\bbR$ is unimodal if there exists $t \in \bbR$ such that $\vp$ is non-decreasing on $\cU \cap (-\infty,t]$ and non-increasing on $[t,+\infty)\cap\cU$.

A symmetric matrix $M$ is said \emph{Robinson} if its rows (equivalently, its columns), are unimodal with respect to their diagonal index. Let $\cR_n$ be the set of Robinson matrix of size $n \times n$. For a vector $\theta = (\theta_0,\dots, \theta_{n-1})\in \bbR^n$, we write $T(\theta)$ for the T\oe plitz matrix with entry $T(\theta)_{i,j} = \theta_{i-j}$, with by convention $\theta_k := \theta_{-k}$ for $k \leq 0$. The matrix $T(\theta)$ is Robinson as soon as the vector $\theta$ is non-increasing. \newclem{For a matrix $X \in \bbR^{n \times n}$ and $i \in [n]$,  $X_i \in \bbR^n$ stands the $i$-th row of $X$.} 

For a permutation matrix $\Pi \in \cS_n$ and any square matrix $M$ of size $n$, we define $\Pi \cdot M := \Pi M \Pi^{\top}$, which corresponds to the action of permuting both columns and rows of $M$ with $\Pi$. In the following, we identify permutation matrices and their corresponding permutations of $[n]$. 


We denote by $C,C',C'',c,c',c'',\dots$ generic numerical constants whose values can differ from one line to another. We will write $x \preceq y$ (resp. $x \succeq y$) for $x \leq Cy$ (resp. $y \geq Cx$). Likewise, we will write $x \propto y$ for $x = Cy$.

\section{Statistical models}\label{sec:models}
We observe a square matrix $Y$ of size $n$ of the form
\beq \label{eq:model}
Y = X + E,
\eeq
where $X$ is a pre-Robinson matrix with entries in $[0,A]$ for some $A \geq 1$. The matrix $E$ is a noise matrix whose entries are independent centered subgaussian variables of subgaussian norm less than $1$. Note that this model covers in particular the case where $Y$ is a matrix of independent Bernouilli random variables $Y_{i,j} \sim \mathrm{Ber}(p_{i,j})$ when the underlying matrix of parameters $(p_{i,j}) \in [0,1]^{n\times n}$ is pre-Robinson.
The goal is to recover a \emph{seriation} of $X$, meaning an ordering $\Pi$ of $[n]$ such that $\Pi \cdot X$ is as close as possible to a Robinson matrix. In this paper, we investigate the recovery in $\ell_2$, leading to the following loss
$$
\ell_2(\Pi) := \inf_{R \in \cR_n} \|\Pi\cdot X - R\|_F,
$$
where $\|\cdot\|_F$ is the Frobenius norm on $n \times n$-matrices. 
Recall that $X$ being pre-Robinson means that there exists a permutation $\Pi^*$ such that $\Pi^* \cdot X$ is Robinson. However, such a permutation is not unique, and this lack of identifiability prevents us from defining the loss as a mere comparison between $\Pi$ and $\Pi^*$. In particular, if $\tau=(n,n-1,\ldots, 2,1)$ is the permutation that reverses the ordering on $[n]$, then both $\Pi^* \cdot X$ and $(\tau \Pi^*)\cdot X$ are Robinson matrices. More generally, the set of permutations that make $X$ a Robinson can be exponentially large in the size of $X$ and its structure can be encoded in a tree, see \cite{atkins1998spectral}. In our context, the lack of idendifiability does not impact our loss --- observe in particular that $\ell_2(\tau\Pi) = \ell_2(\Pi)$ for any permutation $\Pi$.

In this work, we will consider two classes of Robinson matrices:

\paragraph{T\oe plitz model} In this model, we assume that
\beq \label{eq:toepmodel}
X = \Pi^* \cdot T(\theta) \quad \text{where}\quad 
\begin{cases} 
\Pi^* \in \cS_n \quad &\text{and}
\\
\theta \in [0,A]^n \quad &\text{is a nonincreasing vector.} 
\end{cases}
\tag{T\oe}
\eeq
In this context, the loss $\ell_2(\cdot)$ is upper-bounded by the quantity
$$
\ell_2(\Pi) \leq  \|\Pi\cdot X - T(\theta)\|_F.
$$
The T\oe plitz condition ensures some sort of stationarity in the interactions, in the sense that $T(\theta)_{i,j}= \theta_{i-j}$ only depends on the distance $|i-j|$. In fact, T\oe plitz matrices have received a lot of attention in statistics as covariance matrices of stationary time series are T\oe plitz. In the seriation context, those models are in particular used in temporal ordering of single cells in genomics~\cite{karin2023scprisma} or more generally in genome assembly~\cite{cai2023matrix}.

\paragraph{Latent position model}
Let $\vp : \bbR \to [0,A]$ be a symmetric and unimodal function.
In this model, we assume that
\beq \label{eq:latmodel}
X_{i,j} := \vp(V_i-V_j) \quad \text{where}\quad V_1,\dots,V_n \sim \Unif[0,1]\quad \text{iid.} \tag{Lat}
\eeq 
In this context, $X$ is of the form $\Pi^* \cdot M$ where $M$ is Robinson and where $\Pi^*$ is the inverse of the permutation that orders the latent position $V_1,\dots,V_n$ in increasing order. 
The loss $\ell_2(\cdot)$ is then upper-bounded by the quantity
$$
\ell_2(\Pi) \leq  \min\{  \|\Pi\cdot X - M\|_F,\|\tau\Pi\cdot X - M\|_F\}.
$$
Such latent position models are specific instances of the general latent models~\cite{hoff2002latent}. The  model~\eqref{eq:latmodel} also encompasses 1-dimensional random geometric graphs~\cite{Gilbert61,penrose2003random,diaz2020learning}. In the latter, $X$ stands for the adjacency of the graph sampled as follows. For each node $i$, a position $V_i$ is sampled, given $V_i$ and $V_j$, the probability that $i$ is connected to $j$ is equal to $\phi(|V_i-V_j|) \in [0,1]$. 



%

\section{Computational methods for seriation} \label{sec:l2seriation}

\subsection{A generic polynomial algorithm: \texttt{PINES}} \label{sec:generic}

\begin{figure}[!ht]
\centering
\includegraphics[scale = 1.2]{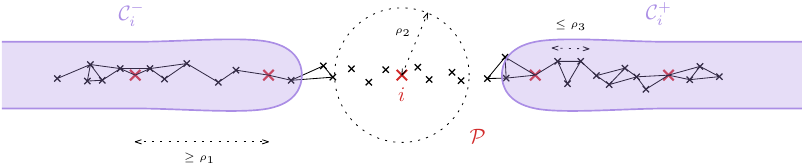}
\caption{A diagram of the construction of the connected components $\cC_i^-$ and $\cC_i^+$. In black crosses are the points of $[n]$ and in red crosses are the points of the $\rho_1$ packing. In black line are the edges of the $\rho_3$ neighborhood graph built on $[n]$ where we removed a $\rho_2$-neighborhood of $i \in \cP$.}
\label{fig:generic}
\end{figure}

We describe in this section an algorithm called \algoname (Polynomial-time Iterative Neighborhood Exploration for Seriation), which intends to seriate a dataset of size $n$ consistently under several seriation models --- among which \eqref{eq:toepmodel} and \eqref{eq:latmodel} --- and which runs in a polynomial time in $n$. \newclem{We assume that there exists an underlying distance on $[n]$, denoted by $\dt$, that is compatible with an optimal ordering in the sense that two indices $i$ and $j$ that are close for this optimal ordering should be close in term of distance $\dt$ --- see \assref{generic} for a more formal description of the requirements. For instance, one could take $\dt(i,j) = |(\Pi^*)^{-1}(i)-(\Pi^*)^{-1}(j)|$ or $\dt(i,j) = \|X_i - X_j\|$. We assume that we have access to an estimator of $\dt$, denoted by $\hat \dt$. We postpone the specific choice of $d$ and $\dt$ to the next subsections.
We will construct an ordering $\Pi$ based on the observations of $\hat \dt(i,j), (i,j) \in [n]$, and assess the accuracy of this seriation in a element-wise fashion with $\dt(\Pi(i),\Pi^*(i))$ for some optimal ordering $\Pi^*$. We refer to $\dt(\Pi(i),\Pi^*(i))$ as \emph{entry-wise} accuracy when considering this seriation procedure. These accuracy measures are then aggregated on the whole dataset --- see \thmref{generic} below.
} The procedure relies on a packing of $[n]$ based on the  empirical distance $\hat \dt$. A $\rho$-packing $\cP$ of $[n]$ for $\hat \dt$ is a subset of $[n]$ such that $\hat \dt(i,j) > \rho$ for all distinct $i,j \in \cP$. A packing is said to be maximal if no superset of $\cP$ can be a packing with respect to the same radius. As a consequence, if $\cP$ is a maximal $\rho$-packing, one can cover $[n]$ with the balls $\{j\in[n]~|~\hat\dt(i,j) \leq \rho\}$ for $i\in \cP$. This induces a partition 
$$
[n] := \bigsqcup_{i \in \cP} Q_i~~~\text{where}~~~i \in Q_i \subset \{j\in[n]~|~\hat\dt(i,j) \leq \rho\}.
$$
The procedure is defined as follows: \newclem{set $\rho_1 > 0$ and
\bitem
\item[1.] Take $\cP$ a maximal $\rho_1$-packing of $[n]$ for $\hat \dt$. This induces a partition of $[n]$;
\item[2.] Order the packing $\cP$ --- see the substeps below;
\item[3.] Seriate all the objects so that  the ordering of the partition is respected. So, we order objects arbitrarily within the sets of the partition.
\eitem 
The idea is to take $\rho_1$ sufficiently large so that we can beat the noise induced by $E$. In doing so, we can find an ordering of $\cP$ which, with high probability,  will align exactly  with the oracle ordering $\Pi^*$. Because the seriation at the level of the packing is perfect, one can expect a entry-wise accurace of order $\rho_1$ in the final estimated permutation, which ultimately yields an $\ell_2$ loss of the order of $n\rho_1$. 

We now elaborate on Step 2: the idea is to build a neighborhood graph on the whole dataset $[n]$, again based on $\hat \dt$. Because of the nature of our models, we expect this graph to behave roughly like a noisy version of an interval graph, which in particular should present one giant connected component containing all the points in the packing. By removing a neighborhood of a point $i \in \cP$ in this graph, we could then expect the giant connected component to be split in two, except if the point $i$ is one of the two extremal points of the packing. This allows us to identify these two extremal points. We can then proceed recursively: if $i_1 \in \cP$ is one extremal point, then the next point of $\cP$ in the ordering should be the such that one of the two connected subcomponent contains no point from $\cP$ but $i_1$. To define this process, we need rely on the scaling parameters $\rho_2$ and $\rho_3$:
\bitem
\item[2.1.] Build a $\rho_3$-neighborhood graph on $[n]$ based on $\hat\dt$;
\item[2.2.] For all $i\in \cP$, consider the subgraph where we removed all points $j \in [n]$ such that $\hat\dt(i,j) \leq \rho_2$;
\item[2.3.] For a suitable choice of parameters, we show that  $\cP \setminus \{i\}$ is contained in at most two connected components of the subgraph. We denote by $\cC_i^-$ and $\cC_i^+$ the trace on $\cP \setminus \{i\}$ of these two connected components (with $\cC_i^+ = \emptyset$ by convention if  $\cP \setminus \{i\}$ is connected in the subgraph). 
\eitem
We then proceed recursively as described above:
\bitem
\item[2.4.] Take a $i_1 \in \cP$ such that $\cC_{i_1}^+ = \emptyset$;
\item[2.5.] Then recursively find $i_k \in \cP$ for $2 \leq k \leq \Card\cP$ such that either $\cC_{i_k}^-$ or $\cC_{i_k}^+$ is $\{i_1,\dots i_{k-1}\}$.
\eitem
The resulting ordering $\Pi_\cP(i_k) := k$ is well defined and a seriation of $\cP$ under the assumptions specified in \assref{generic} and as proven in \thmref{generic}. We refer to \appref{algo} for a more precise pseudo-code writing of the procedure \algoname, and to \figref{generic} for a visual explanation of the procedure. 
}

\begin{ass} \label{ass:generic} We assume that there exists a distance $\dt$ on $[n]$ such that
\bitem
\item[i)] There exists $\ve > 0$ such that
\beq
|\hat \dt(i,j) - \dt(i,j) | \leq \ve,~~~~\forall i,j \in [n].
\eeq
\item[ii)] There exists $\delta > 0$ such that 
\beq
\dt(\Pi^*(i),\Pi^*(i+1)) \leq \delta,~~~~\forall i \in [n-1].
\eeq 
\item[iii)] There exists $\alpha \geq 1$ such that for all $i \leq j \leq k$ in $[n]$
\beq \label{eq:aalpha}
\dt(\Pi^*(i),\Pi^*(k)) \geq \frac1{\alpha} \dt(\Pi^*(i),\Pi^*(j)) ,
\eeq
and the same holds for $k \leq j \leq i$.
\eitem
\end{ass}
\newclem{Let us first elaborate on the assumptions above, and let us do so in reverse order. Point (iii) states that the distance $\dt$ is quasi-increasing for the ordering $\Pi^*$. This means that knowing $\dt$ should allow us, in some sense, to recover some information on $\Pi$, and the quality of this information should decrease as $\alpha$ grows large. We will see in the next Section that this holds for the $L^2$ distance between the rows for some either the \eqref{eq:toepmodel} or \eqref{eq:latmodel} models. Point (ii) states that the distance $\dt$ varies \emph{smoothly} along the ordering $\Pi^*$, in the sense that no gap greater than $\delta$ is allowed. This assumption together with Assumption (iii) implies that two points $i,j \in [n]$ such that $\dt(i,j)$ is small should be close to each other in $\Pi^*$ and vice versa. Finally, Point (i) simply states that we were able to estimate $\dt$ with precision $\ve$. Under these three assumptions, we show that \algoname as described above and for a particular choice of tuning parameters $(\rho_1,\rho_2,\rho_3)$ depending on $\alpha, \delta$ and $\ve$ terminates and we quantify its accuracy in terms of these parameters.}

Recall that $\tau=(n,n-1,\ldots, 2,1)$ is the permutation that reverses the ordering on $[n]$.

\begin{thm} \label{thm:generic} Under \assref{generic}, the procedure \algoname  run with parameters 
$$\rho_3 = \delta + \ve,~~~~ \rho_2= \alpha \delta + 2(1+\alpha) \ve~~~~ \text{and}~~~~ \rho_1= \alpha^2 \delta + (2\alpha^2+3\alpha+2)\ve,$$ 
terminates and yields a permutation $\wh\Pi$ such that, for either $\Pi = \wh\Pi^{-1}$ or $\Pi = \wh\Pi^{-1} \tau$, for all $i \in [n]$, we have
$$
\dt(\Pi(i),{\Pi^*}(i)) \leq (2\alpha+1)\rho_1 + 2\alpha \ve.
$$ In particular, for either $\Pi = \wh\Pi^{-1}$ or $\Pi = \wh\Pi^{-1} \tau$, we have
$$
\(\sum_{i=1}^n \dt(\Pi(i),\Pi^*(i))^p\)^{1/p} \preceq n^{1/p}\alpha^3\(\delta + \ve \),
$$
for all $p \in [1,\infty]$.
\end{thm}

The proof of \thmref{generic} can be found in \appref{proofgeneric}.
\newclem{Let us now comment on the choices of $\rho_3$, $\rho_2$, and $\rho_1$. While running \algoname, we typically want the neighborhood graphs to be such that $\Pi^*(i)$ and $\Pi^*(i+1)$ are connected for all $i \in [n-1]$. This suggests taking $\rho_3 = \ve + \delta$ in light of Points (i) and (ii). The radii $\rho_1$ and $\rho_2$ are then roughly taken proportional to $\delta+\ve$, with the coefficients tuned such that each subgraph has  at most two connected components in the packing. These coefficients depend on $\alpha$ through Point (iii).  In the next sections, we apply this theorem to the models \eqref{eq:toepmodel} or \eqref{eq:latmodel}.
}

\subsection{Seriation for T\oe plitz matrices} \label{sec:modtoe}

In this section, we study Model~\eqref{eq:toepmodel}. The target distance will be
$$
\dt(i,j) := \|X_i-X_j\|.
$$
We show in \prpref{monodistl2} and \prpref{lipdist} of \appref{toeplitz} that, for any non-increasing vector $\theta \in [0,A]^n$, and for any set of indices $1 \leq i < j < k \leq n$, there holds
$$
\|T(\theta)_i - T(\theta)_j\| \leq A \sqrt{|j-i|+1}\quad \text{and}\quad  \|T(\theta)_i - T(\theta)_j\| \geq \frac{1}{\sqrt{2}} \|T(\theta)_i - T(\theta)_k\|,
$$
where we recall that $T(\theta)_i \in \bbR^n$ is the $i$-th row of $T(\theta)$, so that the distance $\dt$ satisfies Points (ii) and (iii) of \assref{generic} with $\delta = \sqrt{2}A$ and $\alpha = \sqrt{2}$.
It only remains to find a suitable candidate for $\hat \dt$. 
Because we do not know the variance of the noise entries $E_{i,j}$ for $i,j \in [n]$, we cannot remove the bias from the estimator of the square of $\dt$ given by $\|Y_i - Y_j\|^2$. In order to circumvent this issue, we write
$$
\dt(i,j)^2 = \|X_i\|^2 + \|X_j\|^2 - 2\inner{X_i}{X_j},
$$
and try to find an estimator for each of these terms.
Because $\inner{Y_i}{Y_j}$ is an unbiased estimator of $\inner{X_i}{X_j}$ we only need to find an estimator for $\|X_i\|^2$. For this, we first define for each $i \in [n]$ a neighborhood of the form
$$
\cN_{i} := \{j \in [n]\setminus\{i\}~|~| S_j - S_i | \preceq A+\sqrt{n \log n}\},
$$
where 
$
S_i := \sum_j Y_{i,j}
$ is an unbiased estimator of $\sum_{j} X_{i,j}$ which coincides with $\|X_i\|_1$ as the entries on $X$ are non-negative.
We then pick
\beq \label{eq:ui}
U_i \in \argmax\{\inner{Y_i}{Y_j}~|~j \in \cN_i\}.
\eeq
We can show that $U_i$ is a good approximant of $\|X_i\|^2$ and that our final estimator of the distances 
\beq \label{eq:hatdunknownvar}
\hat \dt(i,j)^2 := U_i + U_j - 2\inner{Y_i}{Y_j},
\eeq
satisfies the following bound. 
\begin{prp}[Informal]\label{prp:hatdunknownvar}With high-probability, it  holds that
$$
|\hat \dt(i,j) - \|X_i-X_j\| | \preceq A + \sqrt{A}\( n\log(n)\)^{1/4},
$$
uniformly for all $i,j \in [n]$.
\end{prp}
We refer to \appref{proofl2seriation} for more rigorous definitions and proof of this section.
We can now apply \thmref{generic} with $\ve \propto  A + \sqrt{A}\(n\log(n)\)^{1/4}$.

\begin{thm} \label{thm:unknownvar} We let $\wh\Pi$ be the output of \algoname as in \thmref{generic} with $\hat \dt$ defined in \eqref{eq:hatdunknownvar} and with $\alpha = \sqrt{2}$, $\delta = \sqrt{2} A$ and $\ve \propto A + \sqrt{A}\(n\log(n)\)^{1/4}$. The estimator $\wh \Pi$ satisfies
$$
\bbE\left[\ell_2(\wh\Pi)\right] \preceq A \sqrt{n}+ \sqrt{A} n^{3/4} \log^{1/4}(n).
$$
\end{thm}

The proof of this result is also in \appref{proofl2seriation}.

\newclem{
\begin{rem} 
\label{rem:estimA} 
The application of \algoname requires the knowledge of $A$ through the tuning of the parameters. As this quantity is sometimes unknown, one can replace it by the $\max_{i,j \in [n]} Y_{i,j} + \sqrt{8 \log n}$, which is an upper bound of $A$ with probability $\propto 1/n$, yielding only a subsequent loss of a polylog term in the final bound.
\end{rem}
}

\begin{rem} \label{rem:toeplitz} Because a T\oe plitz matrix is invariant if we reverse the orders of the rows and columns, \thmref{generic} actually enforces the bound
$$
\bbE\left[\| \wh \Pi \cdot X - T(\theta)\|_F\right] \preceq A \sqrt{n}+ \sqrt{A} n^{3/4} \log^{1/4}(n) .
$$
\end{rem}
\subsection{Seriation with latent positions}\label{sec:modlat}

In this section, we study the model \eqref{eq:latmodel}. In this context $\|X_i-X_j\|$ does not fulfill the point iii) of \assref{generic} and the targeted distance is chosen to be
$$
\dt(i,j) := \int_{0}^1(\vp_{V_i}(v)-\vp_{V_j}(v)^2dv,\quad\text{where}\quad \vp_s : v \mapsto \vp(v-s).
$$
Notice here that the distance is a random variable. This distance satisfies again point ii) and iii) of \assref{generic} with $\delta \propto A \(\log n\)^{1/2}$ and $\alpha = \sqrt{2}$ with high probability as a consequence of \lemref{eventb} and \lemref{lipg} in the appendix. Like in the previous section, we introduce local neighborhoods for each of the point $i \in [n]$,
$$
\cN_{i} := \{j \in [n]\setminus\{i\}~|~| S_j - S_i | \preceq A \sqrt{n \log n}\} .
$$
Again, the estimator
\beq \label{eq:defhatd2}
\hat \dt^2(i,j) := U_i + U_i - 2 \inner{Y_i}{Y_j},
\eeq
with  $U_i \in \argmax\{\inner{Y_i}{Y_j}~|~j \in \cN_i\}$ satisfies the following properties.
\begin{prp}[Informal]\label{prp:estd2}With high probability, there holds
$$
\sup_{i \neq j} |\hat \dt(i,j) - \dt(i,j) | \preceq A \(n \log(n)\)^{1/4} . 
$$
\end{prp}

We are then in position to apply \thmref{generic}.

\begin{thm} \label{thm:latent} We let $\wh\Pi$ be the output of \algoname  as in \thmref{generic} with $\hat \dt$ defined in \eqref{eq:defhatd2} and with $\alpha = \sqrt{2}$, $\delta \propto A (\log n)^{1/2}$ and $\ve \propto A (n\log n)^{1/4}$. The estimator $\wh \Pi$ satisfies
$$
\bbE\left[ \ell_2(\wh \Pi)\right] \preceq An^{3/4}\log^{1/4}(n).
$$
\end{thm}

\newclem{Let us emphasize that, in this setting, the matrix $X$ is a random matrix, with the randomness stemming from the latent position $V_1,\dots,V_n$. The idea behind the proof is to find deterministic conditions on these latent positions that hold with high-probability and such that, under this condition, Assumption \ref{ass:generic} is met. These conditions boil down to a spacing condition --- see Lemma \ref{lem:eventb} --- and an empirical process result --- see Lemma \ref{lem:eventcd}.}

\section{Computational lower bound} \label{sec:lower:bounds:low_degee}

In this section, we establish a computation-information lower bound of the seriation model with latent positions  using a low-degree polynomial approach --- see~\cite{KuniskyWeinBandeira,WeinSchramm}. This lower bound establishes the optimality of \algoname. We also deal with a variant of the Robinson T\oe plitz model at the end of the section. \\

Let $\lambda\in (0,1)$ be a positive quantity and $k$ be a positive integer that will be fixed later.
We consider a latent position model~\eqref{eq:latmodel} with a function $\phi_{\lambda}$ defined by 
$\phi_{\lambda}(x):= \lambda \mathbf{1}_{|x|\leq 1/\sqrt{n}}$. As in Section~\ref{sec:models}, we can write down our observation model in the form 
\begin{equation}\label{eq:model:CLB}
Y = X+ E, 
\end{equation}
where the entries of $E$ are independent Gaussian variables $\mathcal{N}(0,1)$. Besides, there exists  $\Pi^*$ sampled uniformly on $\mathcal{S}_n$ such that $\Pi^* \cdot X = M$ is some Robinson matrix taking the values $0$ and $\lambda$.
First, we reduce the problem of estimating $\Pi^*$ to that of estimating the matrix $X$. In the following lemma, we show that,  given an estimator $\hat{\Pi}$ of $\Pi^*$, we are able to easily build an estimator of $\hat{X}$ of $X$ with a controlled error. 

\begin{lem}[Reduction to matrix estimation]\label{lem:reduction}
Consider any $\lambda>0$, integer $k$ and any estimator $\widehat{\Pi}$. Define the matrix $\hat{X}$ by $\hat{X}_{i,j}= \lambda/2$ if $|\hat{\pi}^{-1}(i)-\hat{\pi}^{-1}(j)|\leq 2k$ and $\hat{X}_{i,j}= 0$ otherwise. Then, we have 
\begin{equation}\label{eq:upper_bound_risk }
\mathbb{E}\left[\|\hat{X}-X\|_F^2\right] \leq \lambda^2k n + 4 \mathbb{E}[\ell^2_2(\hat{\Pi})] + \lambda^2 n\sqrt{2k}.
\end{equation}
\end{lem}

Assume henceforth that $k\geq 32$. If all polynomial-time estimators $\hat{X}$ satisfy $\mathbb{E}\left[\|\hat{X}-X\|_F^2\right]\geq 1.5 \lambda^2 k n$, then Lemma~\ref{lem:reduction} implies that all polynomial-time estimators of $\Pi^*$ satisfy  $\mathbb{E}[\ell^2_2(\hat{\Pi})]\geq \lambda^2 kn/16$. Hence, it suffices to consider the reconstruction problem of the matrix $X$. As alluded above, we consider the low-degree polynomial framework and we we will establish that no such procedure is able to reconstruct efficiently the matrix $X$.  Given an integer $D>0$, we define $\mathrm{MMSE}_{\leq D}$ as the infimum expected risk achieved by a polynomial estimator of degree up to $D$.
\[
\mathrm{MMSE}_{\leq D}= \inf_{f:\ \mathrm{deg}(f)\leq D}\E[\|f(Y)- X\|^2_{F}]. 
\]
The next theorem states that as long as $\lambda_0$ is small enough and $k$ is small compared to $\sqrt{n}$, no low-degree polynomial estimator of $X$ achieves a small error.

\begin{thm}\label{prp:low_degree_MMSE}
 Define $r_0 = 2\lambda^2(D+1)^4$. If $r_0 <1$ and $(2k+1)^2 \leq n/2$, we have 
    \[
    \mathrm{MMSE}_{\leq D} \geq \lambda^2 2 k (n-1) - 4\lambda^2 k^2 \left(3 + \frac{5r_0}{1-r_0}\right).  
    \]
\end{thm}

Let $\eta>0$ be a positive integer. Let us choose $k$ as the largest integer such that $(2k+1)^2 \leq n/2$ and let us fix $\lambda= 1/(4(\log^{1+\eta}(n)+1)^2)$.  This theorem ensures that,  low-degree polynomials with degree $D\leq (\log(n))^{1+\eta}$ achieve a risk which, up to logarithmic terms, is higher than $n^{3/2}$.
Since lower-bounds for low-degree polynomials with degree $D\leq (\log(n))^{1+\eta}$ are considered~\citep{KuniskyWeinBandeira,WeinSchramm} as evidence of the computational hardness of the problem, Theorem \ref{prp:low_degree_MMSE} suggests computational hardness of estimating $X$ with square Frobenius risk larger than $n^{3/2}$. Then, as a consequence of Lemma~\ref{lem:reduction}, this also suggests the computational hardness of estimating $\Pi^*$ in risk $\ell^2_2$ with a rate no larger than $n^{3/2}$. Since, for any $\Pi$, the loss $\ell_2(\Pi)$ satisfies $\ell_2(\Pi)\leq \|{\bf X}\|_F$, it follows that, with probability higher than $1-1/n^2$, $\ell_2(\Pi)\prec \lambda(kn)^{1/2}$. Recall that we focus on the case where $k$ is of the order of $n^{1/2}$. This implies that an estimator $\hat{\Pi}$ whose $\ell_2^2$ risk is at least of the order of $n^{3/2}$ also has a $\ell_2$ risk at least of the order of $n^{3/4}$. In summary, we have provided evidence, in the low degree computational framework, that no polynomial-time estimator achieves a $\ell_2$ risk significantly faster than $n^{3/4}$ in the latent position model~\eqref{eq:latmodel}. This suggests that the rate $n^{3/4}$ achieved by \algoname{} is optimal among polynomial-time algorithms.

The above low-degree polynomial lower bound has been shown for a latent position model~\eqref{eq:latmodel}. Unfortunately, we are not able to extend Theorem~\ref{prp:low_degree_MMSE} to the Robinson T\oe plitz model~\eqref{eq:toepmodel} for teohnical reasons. Indeed, the proof of Theorem~\ref{prp:low_degree_MMSE} relies on delicate controls of cumulants whose simplifications relies on independences between some of the entries of the matrix $X$; however, such independences do not hold in the  T\oe plitz model. Nevertheless, we are able to show a counterpart of Theorem~\ref{prp:low_degree_MMSE} in a close model defined as follows. Given  a positive integer $k$  and $\lambda\in (0,1)$, define the Robinson-T\oe plitz matrix $M\in \mathbb{R}^{n\times n}$ by $M_{i,j}= \lambda$ if $|i-j|\leq k$ and $M_{i,j}=0$, otherwise. Write $(e_1,\ldots, e_n)$ for the  canonical basis on $\mathbb{R}^n$. Then, we define $\mathcal{S}^{\dag}_n$ as the collection of matrices $\Pi'$  such that each row of $\Pi$ is an element of the canonical basis. In contrast to a permutation matrix, $\Pi'$ can contain identical rows. Such a matrix $\Pi'$ encodes a function $\pi:[n]\mapsto [n]$. Then, we consider the observation model $Y= X + E= \Pi'\cdot M+ E$ where $\Pi'$ is sampled uniformly at random from $\mathcal{S}^{\dag}_n$. For any such $\Pi'$, there exists a permutation matrix $\Pi^*$ such that $\Pi^*\cdot X$ is a Robinson matrix and is, with high probability, close to a T\oe plitz matrix. It is quite straightforward to extend both Lemma~\ref{lem:reduction} and Theorem~\ref{prp:low_degree_MMSE} to this new model and we leave it to the reader. In summary, while we are not able to establish the hardness results for the T\oe plitz model, we give evidence in a slight variation of this model that it is not possible to estimate $\Pi^*$ at the $\ell_2$ rate significantly faster than $n^{3/4}$, which matches the bound in Theorem~\ref{thm:unknownvar}.

\section{Information-theoretic bounds}\label{sec:noncomp}

\subsection{Information-theoretic upper-bounds}
In this section, we provide evidence for the computation-statistical gaps by establishing that the optimal convergence rate in $\ell_2$ distance is of the order of $\sqrt{n}$. For that purpose, we first study a least-square type estimator based on  optimization of criteria over the space of permutations. It is therefore unclear --- and unlikely in the worst case --- that it is possible to efficiently compute them. Nevertheless, we provide them as a benchmark.  
As our main aim in this section is to  show the existence of this gap, we restrict here our attention to the emblematic case where the $E_{i,j}$ are independent Gaussian  variable $\mathcal{N}(0,1)$.
Write $\cA = [0,A] \cap (u\mathbb N)$ for the regular grid of $[0,A]^n$, with grid step $u$.

\paragraph{T\oe plitz model.}
In this paragraph, we consider Model~\eqref{eq:toepmodel}. Let us write $u = 1/n^2$. We consider a least square estimator $\hat \Pi^{\mathrm{(LS,T)}}$ over the grid, as an argmin over $\Pi \in \cS_n$ of
$$\inf_{\tilde \theta\in \cA^n~\mathrm{non-increasing}}\|\Pi \cdot Y - T(\tilde \theta)\|_F^2.$$

\begin{thm}\label{thm:noncompl2}
It holds that
$$\mathbb E[\ell_2(\hat \Pi^{\mathrm{(LS,T)}})]  \preceq \sqrt{n\log\big(n A \big)} + A/n .$$
\end{thm}
The proof of this theorem is in Appendix~\ref{proof:noncomp}.
\paragraph{Seriation with latent positions.} We now consider Model~\eqref{eq:latmodel}.
Set $\mathcal V = [0,1]\cap (\mathbb N/n^8)$ and take $u = 1/n^8$ in the definition of $\mathcal{A}$. We consider a least square estimator $\hat \Pi^{\mathrm{(LS,L)}}$ over a specific grid, as an argmin over $\Pi \in \cS_n$ of
$$\inf_{\tilde V \in \mathcal V^n~\mathrm{non-decreasing},~~\phi \in \mathcal A^n~\mathrm{non-increasing}} \|\Pi \cdot Y - R(\tilde V,\phi)\|_F^2,$$
where $R(\tilde V,\phi)$ is a $n\times n$ matrix such that $R(\tilde V,\phi)_{i,j} = \phi_{|\tilde V_i - \tilde V_j|\times n^8+1}$.

\begin{thm}\label{thm:noncompl2la}
For any vector $v \in [0,1]^n$, it holds that
$$ \mathbb E [\ell_2(\hat \Pi^{\mathrm{(LS,L)}})|V=v]  \preceq \sqrt{n\log\big(n A \big)} + A/n.$$
\end{thm}
The proof of this theorem is in Appendix~\ref{proof:noncomp}. Overall, both Theorems~\ref{thm:noncompl2} and \ref{thm:noncompl2la} imply that, if we set aside computational constraints, it is possible to achieve a $\ell_2$ risk of the order of $\sqrt{n}$.

\subsection{Information-theoretic lower bounds}\label{ss:infLB}

Conversely, we show in this subsection the optimality of the risk $\sqrt{n}$. We again restrict our attention to the case where the $E_{i,j}$ are independent Gaussian  variables $\mathcal{N}(0,1)$.

\paragraph{T\oe plitz model.} In what follows, consider $\theta = u\times  (1,1,0,\ldots,0)$, for some $u>0$ that we will specify later, and write $\mathbb E_\Pi$ for the expectation when the underlying matrix $X$ is $\Pi^\top \cdot T(\theta)$ and write also $\ell_2(\hat \Pi, \Pi):=\ell_2(\hat \Pi)$ (to insist on the dependence of the loss on the true permutation $\Pi$). We have in this context the following lower bound.
\begin{thm}\label{thm:LBl2}
If $u \preceq 1$ and $n \succeq 1$, we have for any estimator $\hat \Pi\in \cS_n$
$$\max_{\Pi \in \mathcal S_n} \mathbb E_\Pi [\ell_2(\hat \Pi)] \succeq u\sqrt{n} .$$  
\end{thm}
This theorem is proven in Appendix~\ref{proof:noncomp}.

\paragraph{Seriation with latent positions.} We now consider Model~\eqref{eq:latmodel}. We have in this context the following lower bound. We consider $\phi = u\times \mathbf 1\{|x|\leq 1/n\}$, for some $u>0$ that we will specify later, and write $\mathbb E$ for the expectation in the associated latent model --- note that $\ell_2(\hat \Pi)$ depends in this case on the latent variables $V$.

\begin{thm}\label{thm:LBl2l}
If $u \preceq 1$, then for any $n_0 \succeq 1$, there exists $n\in [n_0/2, 2n_0]$ such that for any estimator $\hat \Pi\in \cS_n$
$$ \mathbb E [\ell_2(\hat \Pi)] \succeq u\sqrt{n} .$$  
\end{thm}
This theorem is proven in Appendix~\ref{proof:noncomp}.

\section{Discussion} \label{sec:discussion}

\subsection{Seriation with missing values}

A relevant question is whether seriation is still possible in the case of missing values. In this setting, we only observe a mask matrix $B\in \{0,1\}^{n\times n}$ and
$$
Y = B \odot (X + E),
$$
where $X$ is again a shuffled Robinson matrix and the entries of $E$ are iid subgaussian variables of subgaussian norm less than $1$. The symbol $\odot$ denotes the entry-wise multiplication and the entries of $B$ are iid Bernoulli variables of parameter $\lambda \in (0,1]$. If the matrix $B$ is independent from $E$ (in the T\oe plitz model) and from $X$ and $E$, then we have a consistent estimator of the mask parameter
$$
\wh \lambda := \frac1{n^2} \sum_{i,j = 1}^n \ind_{B_{i,j} = 0}.
$$
This estimate allows us to debias the measurements we do from the observation of $Y$. For instance, one could define
\beq \label{eq:hatdmissing}
\hat\dt(i,j)^2 := \frac{1}{\wh\lambda} U_i +  \frac{1}{\wh\lambda} U_j - \frac{2}{\wh\lambda^2} \inner{Y_i}{Y_j},
\eeq 
where $U_i$ would be defined again in the spirit of \eqref{eq:ui} or \eqref{eq:ui2}. In the end, we would get a result of the form.
\begin{thm}[Informal] \label{thm:missing} In the present setting, \algoname  run with $\hat \dt$ defined in \eqref{eq:hatdmissing} and with $\alpha \propto 1$, $\delta \propto A$ and $\ve \propto A {\wh \lambda}^{-3/2} \{n \log n\}^{1/4}$ would output a permutation satisfying
$$
\bbE\left[ \|\wh\Pi\cdot X - M\|_F\right] \preceq  \frac{A}{\lambda^{3/2}} n^{3/4} \log^{1/4}(n),
$$
whenever $\lambda \succeq \sqrt{\log n}/n$, either in \eqref{eq:toepmodel} or \eqref{eq:latmodel}.
\end{thm}


\subsection{Connection to computational barriers in SBM}

The construction of the low-degree polynomial lower bound of Theorem~\ref{prp:low_degree_MMSE} is based on a band matrix $M$ with band size equal to $k$. Then, the signal matrix $X$ is sampled by shuffling the rows and the columns of $X$ according to some function $\pi:[n]\mapsto [n]$ sampled uniformly at random. This construction is reminiscent of the computational barrier~\cite{luo2023computational} that has been recently established for the reconstruction of stochastic block models (SBM) with a large number $K$ of groups. The main difference between our construction and theirs is that, in~\cite{luo2023computational}, the matrix $M$ is block-diagonal with $n/K$ block and the noise is Bernoulli distributed. 

\subsection{Further discussion of existing literature}

\paragraph{Faster Rates under additional assumptions.} We have provided compelling evidence that no polynomial-time estimator  can achieve seriation with a risk much smaller than $n^{3/4}$. However, it is possible to break the computational barrier and to almost achieve the minimax risk under additional conditions either on the affinity function 
$\vp$ or on the non-increasing vector $\theta$. For instance, \cite{issartel2024minimax} considers the case where the vector $\theta$ is bi-Lipschitz which includes the linear case.  In that situation, the authors introduce a polynomial-time estimator achieving the optimal convergence rate $\sqrt{n\log(n)}$. Along those lines, the papers~\cite{janssen2022reconstruction, natik2021consistency} consider some specific Robinson T\oe plitz settings with important spectral gaps and analyze spectral methods under these additional conditions. In a future work, it would be interesting to further characterize the vectors $\theta$ and the affinity functions $\vp$ that allow to bypass the computational barrier $n^{3/4}$. 

\paragraph{Comparison of \algoname and SALB in~\cite{issartel2024minimax}.}
In~\cite{issartel2024minimax}, the authors introduce a procedure SALB that shares some similarities with \algoname. Indeed, the two first steps of SALB amount to first estimating a distance between the rows and using a graph construction similar to 2.2 and 2.3 in Section~\ref{sec:l2seriation} to estimate the set of points which are left or right a given $i\in [n]$. However, there are three important differences between our work and~\cite{issartel2024minimax}. First, we start our procedure by building a maximum packing set on $[n]$, which allows to restrict ourselves to a collection of rows which is possible to seriate. Second, our tuning parameters for \algoname in the T\oe plitz (Theorem~\ref{thm:unknownvar}) and Latent position models (Theorem~\ref{thm:latent}) depend on known quantities such as $A$ and $n$, whereas the tuning parameters in SALB~\cite{issartel2024minimax} depend on the regularity of the matrix $\Pi^*\cdot X$. Finally, we point out that we are able to establish risk bounds for all Robinson T\oe pliz and latent position model whereas~\cite{issartel2024minimax} only consider smooth latent models.

\paragraph{Discussion of Cai and Ma \cite{cai2023matrix}.} The paper~\cite{cai2023matrix} is most related to this work as the authors consider the seriation problem for a permuted T\oe plitz and Robinson matrix, i.e.~Model~\eqref{eq:toepmodel}. However, their objective is different: they investigate separation conditions so that exact reconstruction of the permutation is possible. More precisely, they consider for their model a given subset $\mathcal T_n'$ of the Robinson and T\oe plitz matrices and a given subset $\mathcal S_n'$ of the permutations $\mathcal S_n$, and they define the separation distance of their model as 
$$\rho^*(\mathcal T_n',\mathcal S_n') = \inf_{\Theta \in \mathcal T_n'} \inf_{\Pi_1, \Pi_2 \in \mathcal S_n'} \|\Pi_1\cdot \Theta - \Pi_2\cdot \Theta\|_F.$$
They provide a polynomial-time procedure that achieve exact reconstruction of the permutation with high probability as soon as $\rho^*(\mathcal T_n',\mathcal S_n')$ is at least of the order of $ n^2$ up to logarithmic terms. They also prove information-theoretic upper and lower bounds for their separation rate of the  order of  $ \sqrt{n\log(n)}$.

It is not straightforward to compare their results to ours, as we do not consider the same loss. Nevertheless, we have the following.
\begin{itemize}
    \item If $\rho^*(\mathcal T_n',\mathcal S_n)$ is, up to poly-logarithmic terms, at least of the order of $n^{3/4}$,  then the output of our polynomial-time procedure \algoname achieves perfect recovery with high probability. This is a corollary from our results, \newclem{see Corollary \ref{cor:caima} below}. So that, if no further assumptions are made on the permutations, that is $\mathcal{S}'_n = \mathcal{S}_n / \tau$ (meaning that if $\Pi_1^{-1} \Pi_2 = \tau$ then $\Pi_1$ or $\Pi_2$ is in $\mathcal{S}_n'$ but not both), we significantly outperform  the procedure of~\cite{cai2023matrix} in terms of separation distance for perfect recovery - their rate being of the order of $n^2$.
\newclem{
\begin{cor}\label{cor:caima}
Let $\delta >0$. Consider any given subset $\mathcal T_n'$ of the Robinson and T\oe plitz matrices such that
\begin{equation}\label{eq:sepcai}
    \rho^*(\mathcal T_n',\mathcal S_n) \succeq  \frac1{\delta}(A \sqrt{n}+ \sqrt{A} n^{3/4} \log^{1/4}(n)).
\end{equation} 
If $T(\theta) \in \cT_n'$, then the output of our polynomial-time procedure PINES achieves perfect recovery with probability larger than $1-\delta$.
\end{cor}
\begin{proof}
From Theorem 4 (Remark 5) we have
    $$
\bbE\left[\| \wh \Pi \cdot X - T(\theta)\|_F\right] \preceq A \sqrt{n}+ \sqrt{A} n^{3/4} \log^{1/4}(n).
$$
This implies by Markov inequality that
    $$
\bbP\left(\| \wh \Pi \cdot X - T(\theta)\|_F \succeq \frac1{\delta}(A \sqrt{n}+ \sqrt{A} n^{3/4} \log^{1/4}(n))\right) \leq \delta.
$$
From Equation~\eqref{eq:sepcai}, we know in particular that
$$\inf_{\Pi \in \mathcal S_n \setminus \{\tau, \Id\}} \|\Pi \cdot T(\theta) - T(\theta)\|_F \succeq  \frac1{\delta}(A \sqrt{n}+ \sqrt{A} n^{3/4} \log^{1/4}(n)).$$
Since $X = \Pi^* \cdot T(\theta)$, the two previous equations imply that on an event of probability $1-\delta$:
$$ T(\theta) = \wh \Pi \cdot X = \wh \Pi \Pi^* \cdot T(\theta),$$
so $\ell_2(\wh \Pi) = 0$, namely we perform perfect recovery on this event.
\end{proof}
}
\item Our low-degree lower bound suggests that perfect reconstruction is possible in polynomial time only if  $\rho^*(\mathcal T_n',\mathcal S_n)$  is at least of the order of $n^{3/4}$ --- see Theorem~\ref{prp:low_degree_MMSE} and Lemma~\ref{lem:reduction}. However, we  have no formal proof of this, as we restricted the computational to the related latent model. Nevertheless, we conjecture that the rate $n^{3/4}$ achieved by \algoname{} is optimal for their problem.
\item The tight information-theoretic upper and lower bounds in~\cite{cai2023matrix} for the separation distance are related to our information-theoretic upper and lower bounds in Section~\ref{sec:noncomp} for the Frobenius loss $\ell_2$. However, we only provide results that are tight up to logarithmic factors, unlike~\cite{cai2023matrix} who manage to be tight up to multiplicative constants. 
\end{itemize} 

\paragraph{Discussion of Han et al. \cite{han2023covariance}} The paper~\cite{han2023covariance} considers the problem of covariance alignment, where given two independent samples
$$X_1, \ldots, X_m \sim_{\mathrm{i.i.d.}} \mathcal N(0, \Sigma)~~~~\mathrm{and}~~~~Y_1, \ldots, Y_{m'} \sim_{\mathrm{i.i.d.}} \mathcal N(0, \Pi \cdot \Sigma),$$
where $\Sigma$ is some unknown variance-covariance matrix of dimension $n$ and where $\Pi\in \mathcal S_n$ is some unknown permutation. Their aim is to recover $\Pi$, i.e.~find $\hat \Pi$ such that $\| \Pi \cdot \Sigma -  \hat \Pi \cdot \Sigma\|_F$ is as small as possible. This problem is quite different from ours, however a sub-problem of theirs is related. Indeed, Robinson and T\oe plitz variance-covariance matrices are of special interest as they model well the variance-covariance matrices of some natural stationary processes. In this specific case, the problem in~\cite{han2023covariance} would then be related to our problem, albeit as a two-sample problem and for a different statistical noise structure. Our polynomial-time procedure \algoname could be applied to the estimated variance-covariance matrices constructed based on resp.~$(X_i)_i$ and $(Y_i)_i$, in order to estimate two permutations that, when matched, would give rise to an estimator $\hat \Pi$ of $\Pi$. As the noise structure is different, it is unclear what the error would be, but we conjecture that we would obtain that $\| \Pi \cdot \Sigma -  \hat \Pi \cdot \Sigma\|_F$ is no larger than $n^{3/4}/(m \wedge m')$. This rate is however minimax sub-optimal, as highlighted in~\cite{han2023covariance}. Nevertheless, we believe that our low-degree lower bound might bring some insight on a possible computation-statistical gap in this model which we conjecture to be also of order $n^{3/4}/(m \lor m')$. We leave this interesting open question to future works.

\subsection{Sup-norm seriation under adversarial noise}

The problem of sup-norm seriation consists of finding a permutation $\Pi$ such that
$$
\ell_\infty(\Pi) := \inf_{R \in \cR_n} \| \Pi \cdot X - R\|_\infty,
$$
where $\|\cdot\|_\infty$ is the entrywise sup-norm on the set of matrices, is as small as possible based on the observation of $X+E$ where $X$ is a permuted Robinson matrix and where $E$ is a deterministic (ie. adversarial) noise matrix. Finding $\Pi$ such that $\ell_\infty(\Pi) \leq \|E\|_\infty$ has been shown to be NP-hard by \cite{chepoi2009seriation}. This bounds is unsurprisingly optimal, even in the Toeplitz model, as stated below and proven in \appref{LBinf}

\begin{prp}\label{prp:LBlinf}
There exists a vector $\theta \in \bbR_+^n$ such that, for any estimator $\wh \Pi$ based on $X+E$, it holds that
$$\sup_{\substack{ X = \Pi \cdot T(\theta),~ \Pi \in \cS_n \\ \|E\|_\infty \leq 1}} \|\hat \Pi\cdot X -  T(\theta)\|_\infty = 1.$$
\end{prp}

\begin{figure}[!ht] 
\centering
\includegraphics[scale = .8]{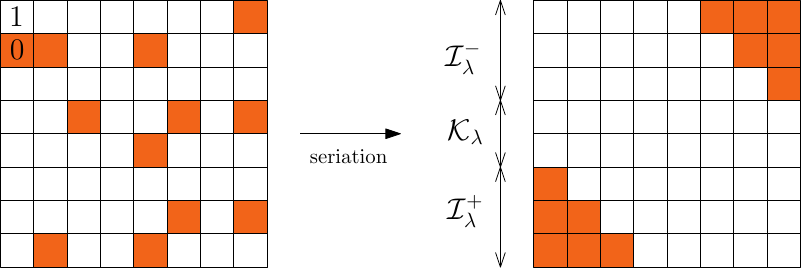}
\caption{(Left) The thresholded matrix $\ind_Q$ before seriation and (Right) the same matrix after seriation, highlighting the partitioning of $[n]$ into three subsets $\cI_\lambda^-$, $\cK_\lambda$ and $\cI_\lambda^+$. The sets $\cI_\lambda^{\pm}$ are already perfectly ordered, and because the increments on $\cK_\lambda$ are bounded, one can apply \algoname  to this subset.}
\label{fig:klambda}
\end{figure}

On the computational side, \cite{chepoi2011seriation} exhibit a polynomial-time algorithm that finds a permutation $\Pi$ whose loss is bounded by $16 \|E\|_\infty$. Although their works takes place in the more general setting of Robinson matrix, we can leverage in our case the algorithm \algoname to get a $c \|E\|_\infty$-approximation of a seriation of $X$ for some numeric constant $c$. An advantage of our method is that it is quite straightforward as compared to the one of \cite{chepoi2011seriation}, although probably yielding a looser approximation (ie. with $c > 16$). In order to apply \algoname, notice that the $\ell_\infty$ loss can be linked to the loss in \thmref{generic} with $\dt(i,j) := \|X_i-X_j\|_{\infty}$ and $p = \infty$. A natural candidate for $\hat\dt$ would then be 
\beq
\label{eq:defhatdlinf}\hat\dt(i,j) := \|Y_i-Y_j\|_{\infty}.
\eeq
It satisfies 
$
| \dt(i,j) - \hat\dt(i,j) | \leq 2\|E\|_\infty.
$
for all $i,j \in [n]$, and it is easy to show that $\dt$ satisfies \eqref{eq:aalpha} for $\alpha = 1$.
However, the increment $\dt(\Pi^*(i),\Pi^*(i+1))$ is not necessarily small (ie of order $\|E\|_\infty$) as required in our analysis of \algoname{} in Section~\ref{sec:l2seriation}. On the other hand, one can expect that big increments of $\dt(\Pi^*(i),\Pi^*(i+1))$ might be helpful to seriate the matrix. We take advantage of this remark by splitting the data in two regions: one where the increment is bounded by a constant (on which we apply \algoname) and one region where the increment is large enough so that the seriation can be determined perfectly. To find this split, we let $\lambda >2\|E\|_\infty$ be a gap parameter and we find the maximal subsets $Q \subset [n]\times [n]$ such that
$$
\min_{(i,j)\in Q} Y_{i,j} \geq \lambda + \max_{(i,j)\in Q^c} Y_{i,j},
$$
and such that $\{Y_{i,j}~|~(i,j)\in Q\}$ contains no gap of size $\lambda$. We then seriate the matrix $\ind_Q$ perfectly using a noiseless algorithm (for instance, \cite{atkins1998spectral}), which splits the data into 3 parts, two of which are already ordered, and the last one being compatible for an application of \algoname. We refer to \figref{klambda} for a diagram of the situation. For a carefully chosen $\lambda$, we are able to derive the following result.

\begin{thm}[Informal] Consider the model $Y=X+E$ where $E$ is deterministic and there exists $\Pi^*$ such that $\Pi^*\cdot X$ is T\oe plitz Robinson. The above algorithm outputs a permutation $\Pi$ such that
$
\ell_\infty(\Pi) \preceq \|E\|_\infty.
$
\end{thm}




\section{Conclusion}

We investigated in this paper the problem of seriating a noisy and permuted Robinson matrix under Frobenius loss. We focused on two natural occurrences of Robinson matrices, namely Toeplitz matrices, and the $1$D latent position model. For both of these instances, we exhibited polynomial-time algorithms with matching lower-bounds on the set of low-degree polynomials. In parallel, we also proved that some non-polynomial time algorithms perform much better, highlighting computational gaps in these settings.

The next natural step would be to provide a polynomial method for seriating general Robinson matrices under $\ell_2$-loss. The fact that polynomial-time algorithms exist in other frameworks such as the aforementioned sup-norm seration \citep{chepoi2011seriation} or in other permutation-estimation problems such as ranking \citep{pilliat2023optimal}, seems to hint towards the existence of such procedures in our context. However, because of the intricate structure of general Robinson matrices, we expect these methods to be much more involved than the one described in this paper.

\section*{Acknowledgements}

The work of N. Verzelen has been partially supported by grant ANR-21-CE23-0035 (ASCAI,ANR). The work of A. Carpentier is partially supported by the Deutsche Forschungsgemeinschaft (DFG)- Project-ID 318763901 - SFB1294 "Data Assimilation", Project A03 and by the DFG on the Forschungsgruppe FOR5381 "Mathematical Statistics in the Information Age - Statistical Efficiency and Computational Tractability", Project TP 02 (Project-ID 460867398). The work of C.~Berenfeld and A.~Carpentier is also partially supported by the Agence Nationale de la Recherche (ANR) and the DFG on the French-German PRCI ANR-DFG ASCAI CA1488/4-1 "Aktive und Batch-Segmentierung, Clustering und Seriation: Grundlagen der KI" (Project-ID 490860858).

\bibliographystyle{abbrv}
\bibliography{ref_bernouli}

\appendix

\section{Probabilistic bounds} \label{app:proba}

This section reviews very basic probability inequalities that we state for sake of completeness and to get explicit constants. 
A centered real-valued random variable $X$ is said to be subgaussian if there exists $\sigma^2$ such that
\beq \label{eq:subg}
 \forall s \in \bbR,~\bbE e^{s X} \leq e^{s^2\sigma^2/2},
\eeq
and we denote $\cS\cG(\sigma^2)$ the set of such random variables. We say that $X$ is sub-exponential \cite[Prp 2.7.1]{vershynin2018high} with parameter $(\sigma^2,\alpha)$ if \eqref{eq:subg} holds but for $|s| \leq 1/\alpha$, and we let $\cS\cE(\sigma^2,\alpha)$ be the set of such random variables. It is straightforward to see that for two independent random variables $X$ and $Y$ in $\cS\cG(\sigma^2)$ (resp. $\cS\cE(\sigma^2,\alpha)$), the sum $X+Y$ is in $\cS\cG(2\sigma^2)$ (resp. $\cS\cE(2\sigma^2,\alpha)$). For the multiplication, we get the following elementary proposition.

\begin{prp} \label{prp:subexpmult} Let $X$ and $Y$ be two independent $\cS\cG(1)$ variables. Then $X^2 - \bbE(X^2) \in \cS\cE(64,4)$ and $XY \in \cS\cE(16,2)$.
\end{prp}
\begin{proof} Using \cite[Thm 2.1]{boucheron2003concentration}, we get first that $\bbE X^2 \leq 4$, and then that for all $|s| \leq 1/4$, letting $\mu = \bbE(X^2)$,
\begin{align*} 
\bbE e^{s (X^2-\mu)} &= 1+\sum_{k=2}^\infty \frac{s^k}{k!} \bbE((X^2-\mu)^{k}) \leq 1+ \sum_{k=2}^\infty \frac{s^k}{k!} (\bbE(X^{2k})+\mu^{k}) \\
&\leq 1+ \sum_{k=2}^\infty \frac{s^k}{k!}  (2 \times 2^k k!+4^k) = 1 + \frac{8 s^2}{1-2s} +(e^{4s}-4s-1) \leq e^{32 s^2}, 
\end{align*} 
so that $X \in \cS\cE(64,4)$. Notice also that the above computation yields
$$
\bbE e^{s X^2} \leq 1+ \frac{8 s^2}{1-2s} \leq e^{16 s^2}.
$$
Then, for any $|s| \leq 1/2$, we have
\begin{align*}
\bbE e^{sXY} \leq \bbE \exp\{\frac{s}2 X^2 + \frac{s}{2} Y^2\} \leq e^{16 (s/2)^2} \times e^{16 (s/2)^2} = e^{8s^2},
\end{align*} 
so that $XY \in \cS\cE(16,2)$, ending the proof.
\end{proof}

\begin{prp} \label{prp:subexpconc} Let $X_1,\dots,X_n$ be $\cS\cE(\sigma^2,\alpha)$ random variables. Then for all $\kappa \geq 0$ such that $\sigma^2 \geq 2\alpha^2(1+\kappa) \log n$, there holds
$$
\bbP\(\max_{1\leq i \leq n} |X_i| \geq \sqrt{2\sigma^2(1+\kappa)\log n}\) \leq 2n^{-\kappa}.
$$
\end{prp}
\begin{proof}Straightforward.
\end{proof}

\section{Proofs of \secref{l2seriation}}\label{app:proofl2seriation}

\subsection{Basic properties of T\oe plitz matrices} \label{app:toeplitz}

We first study the monotonicy of the $L^2$-distance in the T\oe plitz model. Because of boundary effect, we are able to get monotonicity only up to a factor $1/2$.

\begin{prp}[$\ell_2$ monotonicity] \label{prp:monodistl2}
For any indices $i \leq j \leq k$ in $[n]$, there holds
$$
\|M_i - M_k\|^2 \geq \frac12\|M_i - M_j\|^2.
$$
The exact same inequality holds for any indices $k \leq j \leq i$.
\end{prp}



The proof relies on a continuization of $\theta$: we let $\vp : [-1,1] \to \bbR$ be an even function, unimodal with respect to $0$. For any $x\in [0,1]$, we let $\vp_x : s \mapsto \vp(s-x)$ be the version of $\vp$ centered at $x$. It is always well defined on $[0,1]$. We define, for any $x,y \in [0,1]$, the following functions
$$
G_x(t) = \int_0^1 (\vp_x - \vp_t)^2,~~~~F_x(t) = \int_0^t (\vp_x - \vp_t)^2~~~~\text{and}~~~~H_{x,y}(t) = \int_{y}^1 (\vp_x - \vp_t)^2.
$$
We'll show the following properties.

\begin{lem} \label{lem:contprop} For any $x \in [0,1]$, there holds
\bitem
\item[i)] $G_x$ is non-incresing on $[\frac{x}{3},x]$ and non-decreasing on $[x,\frac{x+2}{3}]$;
\item[ii)] $F_x$ is non-decreasing on $[x,1]$;
\item[iii)] $H_{x,y}(t) \geq H_{x,y}\(y\)$ for all $\frac{x+2}{3} \leq y\leq t\leq 1$.
\eitem
\end{lem}

\begin{proof} In the proofs of all three statements, we will assume that $\vp$ is smooth. The results will then follow from the density of smooth function in $L^2([0,1])$. We let $x \in [0,1]$. 

Proof of i): Since $G_x(t) = G_{1-x}(1-t)$, it is sufficient to show that $G_x$ is non-decreasing on $[x,(2+x)/3]$. Since $\vp$ is smooth, so is $G_x$ and for $t > x$, there holds,
\begin{align*}
G_x'(t) &= 2 \int_0^1 \vp'(s-t)(\vp(s-x)-\vp(s-t)) \diff s =:  \int_0^1 g_x(t,s) \diff s \\
&= \int_0^{(x+t)/2} g_x(t,s) \diff s + \int_{(x+t)/2}^{t} g_x(t,s) \diff s + \int_t^{(3t-x)/2} g_x(t,s) \diff s + \int_{(3t-x)/2}^1 g_x(t,s) \diff s
\end{align*}
where we used the fact that $(3t-x)/2 \leq 1$ by assumption. Very simple considerations show that the first and last terms in the last RHS are positive. For the two middle terms, notice that, making the variable change $u = 2t-s$:
\begin{align*}
\int_{(x+t)/2}^{t} g_x(t,s) &=  2 \int_{t}^{(3t-x)/2} \vp'(t-u)(\vp(2t-u-x)-\vp(t-u))) \diff u \\
&= 2 \int_{t}^{(3t-x)/2} \vp'(u-t)(\vp(2t-u-x)-\vp(t-u))) \diff u, 
\end{align*}
so that
$$
\int_{(x+t)/2}^{t} g_x(t,s) \diff s + \int_t^{(3t-x)/2} g_x(t,s) \diff s = 2 \int_{t}^{(3t-x)/2} \vp'(u-t)(\vp(u-x)-\vp(2t-u-x)) \diff u, 
$$
and straight-forward computations show that $|2t-u-x| \leq |u-x|$. 

Proof of ii): For $t \geq x$, there holds
$$
F'_x(t) = 2 \int_0^t \vp'(s-t)(\vp(s-x)-\vp(s-t)) \diff s + (\vp(0)-\vp(s-t))^2.
$$
Now notice that
$$
(\vp(0)-\vp(s-t))^2 = 2\int_{\frac{x+t}{2}}^t (\vp'(s-x)-\vp'(s-t))(\vp(s-x)-\vp(s-t))\diff s
$$
so that
$$
F'_x(t) = 2 \int_0^{\frac{x+t}{2}} \vp'(s-t)(\vp(s-x)-\vp(s-t)) \diff s + 2\int_{\frac{x+t}{2}}^t \vp'(s-x)(\vp(s-x)-\vp(s-t)) \diff s 
$$
and every terms are now positive, hence $F_x$ is non-decreasing on $t \geq x$.

Proof of iii): Let $y \geq \frac{2+x}{3}$. The key point is to notice that $\left|1 - y \right| \leq |x - y|$, so that $\vp_x(u) \leq \vp_t(v)$ for any $t,u,v \in [y,1]$. We now define for $s,t \in [y,1]$ 
$$
\psi_t(s) := \begin{cases}
\vp_t(s)~~~~&\text{if}~~~~ s \geq t \\
\vp_t(1+t-s)~~~~&\text{if}~~~~ s < t 
\end{cases}
$$
which is a rearrangement of $\vp_{y}$ on $[y,1]$. Now since $\vp_x \leq \psi_t \leq \vp_t$ on $[y,1]$, there holds
\begin{align*}
H_{x,y}(t) &=  \int_{y}^1 (\vp_x - \vp_t)^2 \geq  \int_{y}^1 (\vp_x - \psi_t)^2 =  \int_{y}^1 \vp^2_x +  \int_{y}^1 \psi^2_t - 2 \int_{y}^1 \vp_x \psi_t \\
&\geq \int_{y}^1 \vp^2_x +  \int_{y}^1 \vp^2_{y} - 2 \int_{y}^1 \vp_x \vp_{y} = H_{x,y}(y),
\end{align*}
where we used Hardy-Littlewood inequality in the last line.
\end{proof}

\begin{cor} \label{cor:contprop} For any $x \in [0,1]$, and any $t,u \in [0,1]$, there holds
\bitem
\item[i)] If $x \leq t \leq u$, then $G_x(u) \geq \frac12 G_x(t)$;
\item[ii)] If $x \geq t \geq u$, then $G_x(u) \geq \frac12 G_x(t)$.
\eitem
\end{cor}
\begin{proof} Noticing again the symmetry $G_x(t) = G_{1-x}(1-t)$, it is sufficient to show only point $i)$. For $x \leq t \leq u \leq \frac{2+x}{3}$, point i) of \lemref{contprop} applies and the result follows. If $t \leq \frac{2+x}{3} \leq u$, then it is enough to show that $G_x(u) \geq \frac12 G_x(\frac{2+x}{3})$, so that it only remains to show the case when $\frac{2+x}{3} \leq t \leq u$. If $F_x(t) \geq \frac12 G_x(t)$, then, point ii) of \lemref{contprop} yields that $G_x(u) \geq F_x(u) \geq F_x(t) \geq \frac12 G_x(t)$. Otherwise, we have $H_{x,t}(t) = G_x(t) -F_x(t) \geq \frac12 G_x(t)$. But then point iii) of \lemref{contprop} implies that $G_x(u) \geq H_{x,t}(u) \geq H_{x,t}(t) \geq \frac12 G_x(t)$, ending the proof.

\end{proof}


\begin{proof}[Proof of \prpref{monodistl2}] We let $\delta = 1/(2n-1)$, $x_i = 2i\delta$ and define the piecewise constant, even function $\vp : [-1,1] \to \bbR$ with $\vp(x) = \theta_i$ if $x \in (x_i - \delta, x_i+\delta]$ for any $x \geq 0$. 
Now notice that 
\begin{align}
G_{x_i}(x_j) := \int_0^1 (\vp_{x_j}-\vp_{x_i})^2 &= \int_0^{\delta} (\vp_{x_j}-\vp_{x_i})^2 + \sum_{k=1}^{n-1} \int_{x_k-\delta}^{x_k+\delta}(\vp_{x_j}-\vp_{x_i})^2 \nonumber \\
&= 2\delta \|M_i - M_j\|^2- \delta (\theta_i-\theta_j)^2. \label{eq:mphi}
\end{align}
Using that $G_{x_i}(x_k) \geq \frac12 G_{x_i}(x_k)$ thanks to \corref{contprop}, we get that
\begin{align*}
 2\delta \|M_i - M_k\|^2- \delta (\theta_i-\theta_k)^2 \geq \frac{1}{2}\{ 2\delta \|M_i - M_j\|^2- \delta (\theta_i-\theta_j)^2\},
\end{align*} 
so that
\begin{align*} 
\|M_i - M_k\|^2 \geq \frac{1}{2}\|M_i - M_j\|^2+\frac{1}{2}(\theta_i-\theta_k)^2-\frac14(\theta_i-\theta_j)^2 \geq \frac{1}{2}\|M_i - M_j\|^2,
\end{align*} 
we we used that $(\theta_i-\theta_j)^2 \leq (\theta_i-\theta_k)^2$ because $i \leq j \leq k$.
\end{proof}

We end this section with an easy result on a \emph{Lipschitz} property of the $\ell_2$ distance in a T\oe plitz Robinson matrix.

\begin{prp} \label{prp:lipdist}
For any indices $i \leq j$ in $[n]$, if the entries of $\theta$ are in $[0,A]$, there holds
$$
\|M_i - M_j\|^2 \leq A^2 (|i-j| +1).
$$
\end{prp}
\begin{proof} Let $i,j \in [n]$ with $i < j$. Then,
\begin{align*} 
\|M_i - M_j\|^2 &= \sum_{\ell=1}^n (\theta_{\ell-j}-\theta_{\ell-i})^2 \leq A \sum_{\ell=1}^n |\theta_{\ell-j}-\theta_{\ell-i}|\\
&= A \times \{\theta_0-\theta_{j}+\theta_0-\theta_{n-i}+\sum_{\ell = i+1}^{j-1} |\theta_{\ell-j}-\theta_{\ell-i}|\} \\
&\leq A^2(1+|i-j|),
\end{align*} 
ending the proof.
\end{proof}

\subsection{Proofs of \secref{generic}} \label{app:proofgeneric}

For the proof of  \thmref{generic}, we assume without loss of generality that $\Pi^* = \Id$. The proof relies on the following observation.

\begin{prp} \label{prp:connectgeneric} Using the notation of \algoname , and using the parameters of \thmref{generic}, we get, that for all $i \in \cP$, the two sets
$$
V_i^-  = \{j \in \cP~|~j < i\}~~\text{and}~~V_i^+ = \{j \in \cP~|~j > i\},
$$
are contained in two distinct connected components of $G_i$.
\end{prp}
\begin{proof} We first show that $V_i^\pm$ are contained in a connected component. Let for instance $j,k \in V_i^-$ such that $j < k < i$. For any $a [j,k]$, there holds, letting $p = \Pi^*(a)$, $b = \sigma^*(k)$ and $c = \sigma^*(i)$,  that 
\begin{align*} 
\hat \dt(a,i) &\geq \dt(a,i) - \ve \geq \frac{1}{\alpha} \dt(k,i) - \ve \geq \frac12\hat \dt(k,i) - \(1+\frac1{\alpha}\) \ve \geq \frac1{\alpha} \rho_1 -  \(1+\frac{1}{\alpha}\)\ve > \rho_2,
\end{align*} 
so that $a$ is a vertex of $G_i$. Furthermore, if $a \leq k-1$ and, there holds,
$$
\hat\dt(p,p+1) \leq \dt(p,p+1) + \ve \leq \delta + \ve = \rho_3,
$$
so that $(p,p+1)$ is an edge in $G_i$, meaning that $j$ and $k$ are connected.

Last we show that no point of of $V_i^-$ is connected to a point of $V_i^+$. Assume the opposite and find $j \in V_i^-$ and $k \in V_i^+$ that are connected in $G_i$ and let $\gamma = (q_1,\dots,q_L)$ be a path from $j$ to $k$. We let
$$
\ell^- = \sup\{\ell \in [L]~|~ q_\ell < i \}.
$$
and denote $p^- = q_{\ell^-}$, $p^+ = q_{\ell^-+1}$. There holds that $p^- < i < p^+$ and thus
\begin{align*}
\hat\dt(p^-,i) &\leq \dt(p^-,i) + \ve \leq \alpha \dt(p^-,p^+) + \ve \leq \alpha \hat\dt(p^-,p^+) + (1+\alpha)\ve \leq \alpha\rho_3+(1+\alpha)\ve < \rho_2,
\end{align*} 
which is absurd, because $p^-$ is a vertex of $G_i$ so that $\hat\dt(p^-,i) \geq \rho_2$.
\end{proof}

\begin{proof}[Proof of \thmref{generic}] Thanks to \prpref{connectgeneric}, we get that \algoname  terminates and that $\Pi_\cP$ orders $\cP$ according to either $\Id$ of $\tau$. Let assume WLOG that it follows the order of $\Id$. 
We let $a_i^-$ (resp. $a_i^+$) be the minimal (resp. maximal) entry of $P_i$.
Let us first notice that for $j \in \Conv P_i = [a_i^-,a_i^+]$, there holds, if $j \geq i$, $\hat \dt(j,i) \leq \alpha \dt(a_i^+,i) \leq \alpha \hat \dt(a_i^+,i) + \alpha \ve \leq \alpha \rho_1+ \alpha \ve$, and the same inequality but with $a_i^-$ holds if $j \leq i$.

Now let $k \in [n]$ and let $i \in \cP$ such that $\Pi^{-1}(k) \in P_i$. If $k \in \Conv P_i$, then $\dt(k,\Pi^{-1}(k)) \leq  \dt(k,i) + \dt(i,\Pi^{-1}(k)) \leq (\alpha+1) \rho_1+ \alpha \ve$. If $k \notin \Conv P_i$, that means that either all elements of $P_i$ are below or above $k$. If there are below, by cardinality, that means that there exists at least one element $\ell \in [n]$ such that $\ell \in P_j$ with $j < i$ and $\ell > k$. But then there holds that $j < i < k < \ell$ and both $i$ and $k$ are in $\Conv P_j$ so that $\dt(i,k) \leq \dt(i,j) + \dt(j,k) \leq 2 \alpha \rho_1+ 2 \alpha \ve$ and finally $\dt(k,\Pi^{-1}(k)) \leq \dt(k,i) + \dt(i,\Pi^{-1}(k)) \leq (2\alpha+1) \rho_1+ 2\alpha \ve$. The case where all elements of $P_i$ are above $k$ is treated similarly. 


\end{proof}

\subsection{Proofs of \secref{modtoe}}

Again in this section we assume WLOG and for the sake of simplicity that $\Pi^* = \Id$ and thus $M = X$. We let
$$
S_i := \sum_{j=1}^n Y_{i,j} = \sum_{j=1}^n (X_{i,j} + E_{i,j}) = \|X_i\|_1 + R_{i},
$$
where $R_i := \sum_{j=1}^n E_{i,j} \sim \cS\cG(n)$. Furthermore, for $i\neq j$,
$$
\inner{Y_i}{Y_j} =: \inner{X_i}{X_j}+ Q_{i,j}
$$
where $Q_{i,j} = \inner{X_i}{E_j} + \inner{X_j}{E_i} + \inner{E_i}{E_j} \sim \cS\cE(CA^2n,2)$ according to \prpref{subexpmult}. Using \prpref{subexpconc}, for all $\kappa > 0$ we can find event $\cA(\kappa)$ of probability at least $1-4n^{-\kappa}$ such that, as soon as $n \succeq \kappa \log(n)$, 
\begin{align} 
\begin{cases} \displaystyle
\sup_{1 \leq i \leq n} |R_i | &\preceq \sqrt{\kappa n\log(n)} =: \nu_n(\cA_1),  \\
\displaystyle
\sup_{1 \leq i \neq j \leq n} |Q_{i,j}| &\preceq A\sqrt{\kappa n\log(n)} =: \nu_n(\cA_2),
\end{cases} \label{eq:rq}
\end{align} 
For any pair $i,j \in [n]$, assuming that $\|X_j\|_2 \geq \|X_i\|_2$, there exists a set of indices $\cL$ such that
$$
\|X_j\|^2-\|X_i\|^2 = \sum_{\ell \in \cL} \theta_\ell^2 - \theta_{\ell+|i-j|}^2 
~~~~\text{and}~~~~
\|X_j\|_1-\|X_i\|_1 = \sum_{\ell \in \cL} \theta_\ell - \theta_{\ell+|i-j|},
$$
so that in particular $\|X_j\|^2-\|X_i\|^2 \leq 2A(\|X_j\|_1-\|X_i\|_1)$. For $i \in [n]$, we define 
$$
\cN_{i} := \{j \in [n]\setminus\{i\}~|~| S_j - S_i | \leq 2A+2 \nu_n(\cA_1)\}.
$$
This set is a subset of 
$$
 \{j \in [n]\setminus\{i\}~|~| \|X_i\|^2 - \|X_j\|^2 | \leq 4A^2+8A \nu_n(\cA_1)\}
$$
and contains at least $k = i-1$ or $k = i+1$ because for all $k \in [n-1]$, there holds $|\|M_{k}\|_1-\|M_{k+1}\|_1 | \leq 2A$. Finally, let
\beq \label{eq:uisupp} \notag
U_i := \argmax\{\inner{Y_i}{Y_j}~|~j \in \cN_i\}.
\eeq
We show that $U_i$ is a good approximant of $\|X_i\|^2$.

\begin{prp} \label{prp:ui} On the event $\cA(\kappa)$, there holds
$$
| U_i - \|X_i\|^2 | \preceq A^2 + A \sqrt{\kappa n\log(n)},
$$
uniformly for all $i \in [n]$.
\end{prp}
\begin{proof} We let $k \in [n]$ such that $i = \Pi^*(k)$ and assume WLOG that $k < n$. We have $j = \Pi^*(k+1)$ in $\cN_i$ and $|\|M_{k}\|^2-\|M_{k+1}\|^2 | \leq 2A^2$. Now notice that on the event $\cA(\kappa)$,
\begin{align*}
U_i &\geq \inner{Y_i}{Y_j} \geq \inner{X_i}{X_j} -  \nu_n(\cA_2) \\
&= \frac12\{\|X_i\|^2+\|X_j\|^2-\|X_i-X_j\|^2\}- \nu_n(\cA_2) \\
&\succeq \|X_i\|^2 - A^2 - \nu_n(\cA_2),
\end{align*} 
where we used \prpref{lipdist}. Furthermore, for all $j \in \cN_i$,
\begin{align*} 
\inner{Y_i}{Y_j} &\leq  \inner{X_i}{X_j} + \nu_n(\cA_2) \\
&\leq \frac{1}{2} \{\|X_i\|^2 + \|X_j\|^2\} + \nu_n(\cA_2) \\
&\leq \|X_i\|^2+ \frac{1}{2}\{4A^2+8A \nu_n(\cA_1)\}+\nu_n(\cA_2) 
\end{align*} 
so that
$$
U_i \preceq \|X_i\|^2 + A^2 + A \nu_n(\cA_1) + \nu_2(\cA_2),
$$
which ends the proof.
\end{proof}

\begin{prp} \label{prp:hatdunknownvarsupp}On the event $\cA(\kappa)$, there holds
$$
|\hat \dt(i,j) - \|X_i-X_j\| | \preceq A + \sqrt{A}\(\kappa n\log(n)\)^{1/4},
$$
uniformly for all $i,j \in [n]$.
\end{prp}

\begin{proof} We get straightforwardly that 
$$|\hat \dt(i,j)^2 - \|X_i-X_j\|^2 | \preceq A^2 + A \sqrt{\kappa n \log(n)},$$ 
for all $i,j \in [n]$. If $\|X_i-X_j\|^2 \preceq A^2 + A \sqrt{\kappa n \log(n)}$, then $\hat \dt(i,j)^2 \preceq A^2 + A \sqrt{\kappa n \log(n)}$ and thus
\begin{align*} 
|\hat \dt(i,j) - \|X_i-X_j\|| \preceq \{A^2 + A \sqrt{\kappa n \log(n)}\}^{1/2} \preceq A + \sqrt{A} \{\kappa n \log(n)\}^{1/4}.
\end{align*} 
Now if $\|X_i-X_j\|^2 \succeq A^2 + A \sqrt{\kappa n \log(n)}$, we find that
\begin{align*} 
|\hat \dt(i,j) - \|X_i-X_j\| | &= \frac{|\hat \dt(i,j)^2 - \|X_i-X_j\|^2 |}{\hat \dt(i,j) + \|X_i-X_j\|} \preceq \{A^2 + A \sqrt{\kappa n \log(n)}\}^{1/2} \\
&\preceq A + \sqrt{A} \{\kappa n \log(n)\}^{1/4},
\end{align*} 
which ends the proof.
\end{proof}

\begin{proof}[Proof of \thmref{unknownvar}] Applying \thmref{generic}, we get that, on the event $\cA(\kappa)$, \algoname  terminates and yields a permutation $\wh \Pi$ such that
$$
\|\wh\Pi \cdot X - X\|_F \preceq  A\sqrt{n} + \sqrt{A} n^{3/4}\{\kappa \log(n)\}^{1/4}.
$$
We then find, taking $\kappa \geq 1/2$:
\begin{align*} 
\bbE \|\wh\Pi \cdot X - X\|_F &\leq A\sqrt{n} + \sqrt{A} n^{3/4} \{\kappa \log(n)\}^{1/4} + 2An \bbP(\cA(\kappa)^c) \\
&\preceq A\sqrt{n} + \sqrt{A} n^{3/4} \{\kappa \log(n)\}^{1/4},
\end{align*} 
which ends the proof.
\end{proof}

\subsection{Proofs of \secref{modlat}} \label{app:proofmodlat}

We recall the notation for $s,t \in [0,1]$,
$$
G_t(s) := \int_{0}^1(\vp_s-\vp_t)^2,\quad\text{where}\quad \vp_s : v \mapsto \vp(v-s).
$$
We let $\dt(i,j)^2 := nG_{V_i}(V_j) = nG_{V_j}(V_i)$. This distance satisfies, thanks to \corref{contprop},
$$
\dt(\Pi^*(i),\Pi^*(k)) \geq \frac{1}{\sqrt{2}} \dt(\Pi^*(i),\Pi^*(j)) \quad \forall i\leq j \leq k.
$$
We let $W_i = V_{{\Pi^*}^{-1}(i)}$ be the ordered latent positions. 
\begin{lem} \label{lem:eventb}
The event $\cB({\kappa})$ for which 
$$\forall i \in [n-1],\quad |W_i - W_{i+1}| \preceq  \frac{\kappa \log n}{n}  =: \nu_n(\cB)$$ 
has probability at least $1-n^{-\kappa}$.     
\end{lem}

\begin{proof}[Proof of \lemref{eventb}]
Let $N > 1$ and $I_k = [(k-1)/N,k/N]$ for $k \in [N]$. There holds
$$
\bbP(\exists k \in [N], I_k \cap \{U_1,\dots,U_n\} = \emptyset) \leq N (1-1/N)^n \leq N e^{-n/N}.
$$
For $N = \floor{n / (\kappa + 1) \log n}$, and $n \geq 2 (\kappa + 1) \log n$, we find that
$$
\bbP(\exists k \in [N], I_k \cap \{U_1,\dots,U_n\} = \emptyset) \leq n^{-\kappa} 
$$
and that on the complementary of this event,
$$
\forall i \in [n-1],\quad |W_i - W_{i+1}| \leq \frac1N \preceq \frac{\kappa \log n}{n}.
$$
\end{proof}
\begin{lem} \label{lem:lipg}For all $s,t \in [0,1]$, 
$$G_t(s) \leq A \int |\vp_t-\vp_s| \leq 2A^2|s-t|.$$
\end{lem}
\begin{proof} Assume WLOG that $t>s$. We simply write, letting $v = (x+t)/2$ and $a = t-s$,
\begin{align*}
G_t(s) &\leq A \int_0^1 |\vp_s -\vp_t| = A \int_{0}^{v} \vp_s - \vp_t + \int_v^1 \vp_t-\vp_s \\
&=A \{ \int_{-a}^{v-a} \vp_t - \int_0^v \vp_t + \int_v^1 \vp_t - \int_{v-a}^{1-a} \vp_t\} \\
&\leq A\{\int_{-a}^0 \vp_t + \int_{1-a}^a \vp_t\} \leq 2A^2 a.
\end{align*}
\end{proof}
On the event $\cB(\kappa)$, the distance $\dt$ satisfies, for all $i \in [n-1]$,
\begin{align*}
\dt(\Pi^*(i),\Pi^*(i+1))^2 &=  n \int_0^1 \(\vp_{W_i} - \vp_{W_{i+1}}\)^2 \leq 2A^2 n|W_i - W_{i+1}| \preceq 2A^2 n \nu_n(\cB).
\end{align*} 
where we used \lemref{lipg}. It only remains to find an estimator of $\dt$. We proceed exactly like in the previous section.
We let
\begin{align*}
\cF &:= \{\vp_t, t\in[0,1]\} \quad \text{and} \quad 
\cG := \{(\vp_s-\vp_t)^2, (s,t) \in [0,1]^2\} \bigcup \{\vp_s \vp_t, (s,t) \in [0,1]^2\}. 
\end{align*}
\begin{lem} \label{lem:eventcd} If $P = \Unif [0,1]$, and for $P_n$ being the empirical measure associated with a $n$-sample of $P$, there holds
\begin{enumerate}
    \item with probability at least $1-n^{-\kappa-1}$,
\beq \label{eq:evf}
n \sup_{f \in \cF} |P_n(f)-P(f)|  \preceq A \sqrt{\kappa n \log(n)} =: \nu_n(\cC);
\eeq
    \item with probability at least $1-n^{-\kappa-2}$,
\beq \label{eq:evg}
n \sup_{f \in \cG} |P_n(f)-P(f)| \preceq A^2 \sqrt{\kappa n \log(n)} =: \nu_n(\cD).
\eeq
\end{enumerate}

\end{lem}

\begin{proof}[Proof of \lemref{eventcd}] We let
$$
Z_1 = n \times \sup_{f \in \cF} |P_n(f)-P(f)| \quad \text{and} \quad Z_2 = n \times \sup_{f \in \cG} |P_n(f)-P(f)|
$$
By Talagrand's inequlity \cite[Thm 2.3]{bousquet2002bennett}, there holds, for $k \in \{1,2\}$ and all $x \geq 0$,
$$
\bbP(Z_k \geq \bbE Z_k + \sqrt{(A^{2k} n + 2A^k \bbE Z_k) x} + A^k x/3) \leq e^{-x}.
$$
Furthermore, we can bracket the set $\{\vp_s, s\in[t,u]\}$ with $t<u$ by $[\underline{\vp},\overline{\vp}]$ where 
$$
\underline{\vp}(x) := \begin{cases}
\vp_u(x) \quad &x < t \\
0 \quad &x \in [t,u]\\
\vp_t(x) \quad &x > u
\end{cases}
\quad
\text{and}
\quad
\overline{\vp}(x) := \begin{cases}
\vp_t(x) \quad &x < t \\
A \quad &x \in [t,u]\\
\vp_u(x) \quad &x > u.
\end{cases}
$$
This bracket is such that
$$
\int (\underline{\vp}(x) - \overline{\vp}(x) )^2 \leq A^2 |t-u|+\int_0^1 (\vp_t-\vp_u)^2 \leq 3A^2|t-u|,
$$
where we used \lemref{lipg}. We deduce that $N_{[]}(\cF,L^2(P),\ve) \leq \max\{3A^2/\ve,1\}$. Likewise, we would get that $N_{[]}(\cG,L^2(P),\ve) \leq \max\{cA^4/\ve^2,1\}$. Using for instance \cite[Thm 3.5.13]{gine2021mathematical}, we find that
$$
\bbE Z_k \preceq A^k \sqrt{n}\quad \forall k\in \{1,2\}.
$$
Applying Talagrand's inequality to $x = \kappa \log n$ thus yields
$$
Z_k \preceq A^k \sqrt{\kappa n \log n},
$$
with probability at least $1-n^{-\kappa}$ for both $k \in \{1,2\}$.
\end{proof}

We let $\cC_{i}(\kappa)$ be the be the corresponding event of \eqref{eq:evf} for $P = \Unif[0,1]$ and
$$P_n^{(i)} := \frac{1}{n-1} \sum_{k \neq i} \delta_{X_k}.$$ 
and $\cD_{i,j}(\kappa)$ be the corresponding event of \eqref{eq:evg} for $P = \Unif[0,1]$ and
$$P_n^{(i,j)} := \frac{1}{n-2} \sum_{k \neq i,j} \delta_{X_k}.$$ 
Like under the model \eqref{eq:toepmodel}, the event $\cA(\kappa)$ still has probability at least $1-4n^{-\kappa}$ (to see that, simply work conditionally to $\{V_1,\dots,V_n\}$).
We set the rest of the analysis on the event 
$$\cE(\kappa) := \cA(\kappa) \cap \cB(\kappa) \cap  \bigcap_{i} \cC_i(\kappa) \cap  \bigcap_{i \neq j} \cD_{i,j}(\kappa),$$
which has probability at least $1- 7n^{-\kappa}$. We can then also write that 
$$
S_i := \sum_{j=1}^n Y_{i,j} = \sum_{j=1}^n X_{i,j}+E_{i,j},
$$
so that, on $\cE(\kappa)$,
$$
\sup_{i}\left|S_i - n \int \vp_{V_i}\right| \leq \nu_n(\cA_1) + \nu_n(\cC). 
$$
Likewise, 
$$
\sup_{i \neq j} \left|\inner{Y_i}{Y_j}- n \int \vp_{V_i} \vp_{V_j}\right| \leq \nu_n(\cA_2) +  \nu_n(\cD). 
$$
Because again for all $s,t \in [0,1]$
$$
\int \vp_s^2 - \int \vp_t^2 \leq A \int \vp_s - A \int \vp_t,
$$
we can still introduce
$$
\cN_{i} := \{j \in [n]\setminus\{i\}~|~| S_j - S_i | \leq 2An\nu_n(\cB)+2\nu_n(\cA_1) + 2 \nu_n(\cC)\},
$$
and
\beq \label{eq:ui2}
U_i \in \argmax \{\inner{Y_i}{Y_j}~|~j \in \cN_i\}.
\eeq
\begin{prp} \label{prp:ui2} On the event $\cE(\kappa)$, there holds
$$
\left| U_i - n \int \vp_{V_i}^2 \right| \preceq A^2 \sqrt{\kappa n \log(n)},
$$
uniformly for all $i \in [n]$.
\end{prp}

\begin{proof}[Proof of \prpref{ui2}] We proceed like in the proof of \prpref{ui}.
We let $k \in [n]$ such that $i = \Pi^*(k)$ and assume WLOG that $k < n$. We have $j = \Pi^*(k+1)$ in $\cN_i$ because
\begin{align*}
|S_i - S_j| &\leq n|P \vp_{V_i} - P \vp_{V_j}| + 2 \nu_n(\cC) + 2 \nu_n(\cA_1) \\
&\leq 2A n\nu_n(\cB) + 2 \nu_n(\cC) + 2 \nu_n(\cA_1),
\end{align*}
where we used \lemref{lipg}. Now notice that
\begin{align*}
U_i &\geq \inner{Y_i}{Y_j} \geq \inner{X_i}{X_j} -  \nu_n(\cA_2) \\
&\geq n \int \vp_{V_i} \vp_{V_j} - \nu_n(\cD) - \nu_n(\cA_2) \\
&=  n \int \vp_{V_i}^2 +  n \int\vp_{V_i}\(\vp_{V_j}-\vp_{V_i}\) - \nu_n(\cD) - \nu_n(\cA_2)
\\
&\succeq  n \int \vp_{V_i}^2 - A^2 n \nu_n(\cB) - \nu_n(\cD) - \nu_n(\cA_2).
\end{align*} 
Furthermore, for all $j \in \cN_i$,
\begin{align*} 
\inner{Y_i}{Y_j} &\leq  \inner{X_i}{X_j} + \nu_n(\cA_2) \leq n \int \vp_{V_i} \vp_{V_j} + \nu_n(\cD) + \nu_n(\cA_2) \\
&= \frac{n}{2} \{\int \vp_{V_i}^2  + \int \vp_{V_j}^2\} + \nu_n(\cD) + \nu_n(\cA_2) \\
&\leq n \int \vp_{V_i}^2+ \frac{A}{2}\{2A n\nu_n(\cB) + 4 \nu_n(\cC) + 4 \nu_n(\cA_1)\}+ \nu_n(\cD) + \nu_n(\cA_2)
\end{align*} 
so that
$$
U_i \preceq n \int \vp_{V_i}^2 + A^2 n\nu_n(\cB) + A \nu_n(\cC) + A \nu_n(\cA_1) + \nu_n(\cD) + \nu_n(\cA_2),
$$
which ends the proof.    
\end{proof}

\begin{prp}\label{prp:estd2supp} On the event $\cE(\kappa)$, there holds
$$
\sup_{i \neq j} \left|\hat \dt(i,j) - \dt(i,j) \right | \preceq A \{\kappa n \log(n)\}^{1/4}.
$$
\end{prp}

\begin{proof}[Proof of \prpref{estd2}] We proceed like in the proof of \prpref{hatdunknownvar}. 
We get straightforwardly that 
$$|\hat \dt(i,j)^2 - \dt(i,j)^2 | \preceq  A^2 \sqrt{\kappa n \log(n)},$$ 
for all $i,j \in [n]$. If $\dt(i,j)^2 \preceq A^2 \sqrt{\kappa n \log(n)}$, then $\hat \dt(i,j)^2 \preceq A^2 \sqrt{\kappa n \log(n)}$ and thus
\begin{align*} 
|\hat \dt(i,j) -\dt(i,j)| \preceq \{A^2 \sqrt{\kappa n \log(n)}\}^{1/2}.
\end{align*} 
Now if $\dt(i,j)^2 \succeq A^2 \sqrt{\kappa n \log(n)}$, we find that
\begin{align*} 
|\hat \dt(i,j) - \dt(i,j) | &= \frac{|\hat \dt(i,j)^2 - \dt(i,j)^2 |}{\hat \dt(i,j) + \dt(i,j)} \preceq \{A^2 \sqrt{\kappa n \log(n)}\}^{1/2},
\end{align*} 
which ends the proof.
\end{proof}


\section{Proofs of \secref{lower:bounds:low_degee}}

\subsection{Proof of Lemma~\ref{lem:reduction}}

If $|\hat{\pi}^{-1}(i)- \hat{\pi}^{-1}(j)|\leq 2k$, we have $|\hat{X}_{i,j}-X_{i,j}|= \lambda/2$. Hence, it follows that 
\begin{equation}\label{eq:upper_bound_l2_pi3}
\|\hat{X}-X\|_F^2\leq \lambda^2 kn + \lambda^2 \left|\{(i,j): X_{i,j}= \lambda \text{ and } |\hat{\pi}^{-1}(i)- \hat{\pi}^{-1}(j)|\geq 2k\}\right| 
\end{equation}
We will bound the RHS using the loss $\ell^2_2(\hat{\Pi})$. For that purpose, let us introduce $\mathcal{R}'_n$ the collection of $n\times n$ matrices whose rows are unimodal and achieve their maximum on the diagonal. Obviously $\mathcal{R}_n\subset \mathcal{R}'_n$. Besides, we introduce $\mathcal{R}'_n(\lambda)$ as the subset of $\mathcal{R}_n$ that only take its values in $\{0,\lambda\}$. We have
\begin{equation}\label{eq:lower_bound_l2_pi}
\ell^2_2(\hat{\Pi})= \inf_{R\in \mathcal{R}_n}\|\hat{\Pi}\cdot X - R \|_F^2 \geq  \frac{1}{4}\inf_{R\in \mathcal{R}'_n(\lambda)}\|\hat{\Pi}\cdot X - R \|_F^2 .
\end{equation}
Let us fix $R$ as any matrix in $\mathcal{R}'_n(\lambda)$ that achieves the above infimum.  Fix any $i\in [n]$. We claim that  
\begin{equation}\label{eq:lower_row}
\|[\hat{\Pi}\cdot X - R]_{i}\|^2\geq \lambda^2 \left|\{j: X_{\hat{\pi}(i),\hat{\pi}(j)}= \lambda \text{ and } |i-j|\geq 2k\}\right| - \lambda^2 \left[\left|\{j: X_{\hat{\pi}(i),\hat{\pi}(j)}= \lambda\}\right|  - (2k+1) \right]_+ ,  
\end{equation}
where $[x]_+= \max(x,0)$. Let us show~\eqref{eq:lower_row}. First, the inequality is trivial if the rhs is non-positive. Next, at worst all the non-zero entries $X_{\hat{\pi}(i),\hat{\pi}(j)}$ are on the same side with respect to $\hat{\pi}(i)$. By symmetry, we assume henceforth that there are on the right side. Also, there are at least 
\[
\left|\{j: X_{\hat{\pi}(i),\hat{\pi}(j)}= \lambda \text{ and } |i-j|\geq 2k\}\right| - \lambda^2 \left[\left|\{j: X_{\hat{\pi}(i),\hat{\pi}(j)}= \lambda\}\right|  - (2k+1) \right]_+\ , 
\]
entries of $(X_{\hat{\pi}(i),\hat{\pi}(i)+1},\ldots X_{\hat{\pi}(i),\hat{\pi}(i)+2k})$ that are equal to zero.  As a consequence, if $R_{\hat{\pi}(i),\hat{\pi}(i)+2k}=\lambda$, then~\eqref{eq:lower_row} holds. Otherwise, if $R_{\hat{\pi}(i),\hat{\pi}(i)+2k}=0$, this implies that $\|[\hat{\Pi}\cdot X - R]_{\hat{\pi}(i)}\|^2/\lambda^2$ is at least the number of non-zero entries of  $(\hat{\Pi}\cdot X)_i$ that are at distance larger than $2k$ from $i$. We have shown~\eqref{eq:lower_row}.
It then follows from~\eqref{eq:lower_row} and~\eqref{eq:lower_bound_l2_pi} that 
\[
4\ell^2_2(\hat{\Pi})\geq \lambda^2 \left|\{(i,j): X_{i,j}= \lambda \text{ and } |\hat{\pi}^{-1}(i)- \hat{\pi}^{-1}(j)|\geq 2k\}\right| - \lambda^2\sum_{i=1}^n\left[\left|\{j: X_{i,j}= \lambda\}\right|  - 2k \right]_+ 
\]
Together with \eqref{eq:upper_bound_l2_pi3}, we obtain that 
\[
\|\hat{X}-X\|_F^2\leq \lambda^2 kn + 4\ell^2_2(\hat{\Pi}) + \lambda^2\sum_{i=1}^n\left[\left|\{j: X_{i,j}= \lambda\}\right|  - 2k \right]_+ \ .  
\]
To conclude the proof, it suffices to control the expectation of the last term. For a fixed $i$, $|\{j: X_{i,j}= \lambda\}|-1$ is stochastically upper bounded by a Binomial distribution with parameters $n$ and $2k/n$. By Jensen inequality, we deduce that $\mathbb{E}\left[\left|\{j: X_{i,j}= \lambda\}\right|  - (2k+1) \right]_+\leq \sqrt{2k}$. The result follows.

\subsection{Proof of Theorem~\ref{prp:low_degree_MMSE}}

First, we reduce the problem of estimating $X$ to that of estimating the functional $x=X_{1,2}$. Since the diagonal of $X$ is almost surely equal to $\lambda$ and since the random variable $(X_{i,j})_{i<j}$ is exchangeable it follows that 
\[
\mathrm{MMSE}_{\leq D}= \inf_{f:\ \mathrm{deg}(f)\leq D}\E[\|f(Y)- X\|^2_{F}] = n(n-1)  \inf_{g:\ \mathrm{deg}(g)\leq D}\E[(g(Y)- x)^2]\ . 
\]
Hence, we focus on lower bounding the $\mathrm{MMSE}'_{\leq D}=  \inf_{g:\ \mathrm{deg}(g)\leq D}\E[(g(Y)- x)^2]$. First, we compute the first moment of $x$, which corresponds to the probability that any two points have their latent label at a distance smaller than $k$. 
\[
    \E[x]= \lambda\left[\frac{2k}{n}- 2\int_0^{k/n}tdt\right] =  \lambda\left[\frac{2k}{n}- \frac{k^2}{n^2}\right] =: \lambda p_0 \ .
\]
Henceforth, we define $\rho = 2k/n$, which is an upper bound of $p_0$.
This proof is based on the general technique of~\cite{WeinSchramm} for lower bounding the $ \mathrm{MMSE}_{\leq D}$ in signal + noise Gaussian model. In particular, it is established in~\cite{WeinSchramm} that 
\begin{equation}\label{eq:decomposition_MMSE}
\mathrm{MMSE}'_{\leq D}= \mathbb{E}[x^2]- \mathrm{Corr}^2_{\leq D}=  \lambda^2p_0 - \mathrm{Corr}^2_{\leq D}    \ , 
\end{equation}
where the low-degree correlation is defined by 
\[
\mathrm{Corr}_{\leq D}= \sup_{\deg f \leq D}\frac{\E_{(x,Y)}[f(Y)x]}{\sqrt{\E_{Y\sim \P}[f(Y)^2] }}\ . 
\]
Since $\mathbb{E}[x^2]=\lambda\mathbb{E}[x]$, we only have to bound $\mathrm{Corr}_{\leq D}$.
Given a matrix $\alpha\in \mathbb{N}^{n\times n}$ with integer values we henceforth write $|\alpha|= \sum_{i,j}\alpha_{i,j}$ for the sum of its entries, $\alpha!= \prod_{i,j} \alpha_{i,j}!$, and $X^{\alpha}= \prod_{i,j}X_{i,j}^{\alpha_{i,j}}$. Besides, for two such matrix $\beta$ and $\alpha$, we write that $\beta\leq \alpha$ if the inequality holds entry wise and $\binom{\alpha}{\beta}= \prod_{i,j} \binom{\alpha_{i,j}}{\beta_{i,j}}$.  By Theorem 2.2 in~\cite{WeinSchramm}, we have 
\begin{eqnarray}\label{eq:upper:low_correlation}
\mathrm{Corr}_{\leq D}^2 \leq \sum_{\alpha \in \mathbb{N}^{n\times n},\, 0\leq |\alpha|\leq D}\frac{\kappa_{\alpha}^2}{\alpha!}\ , 
\end{eqnarray}
where the quantity $\kappa_{\alpha}$ is defined recursively by 
\[
    \kappa_{\alpha}= \mathbb{E}[xX^{\alpha}]- \sum_{0\leq \beta \lneq \alpha }\kappa_{\beta}\binom{\alpha}{\beta}\mathbb{E}[X^{\alpha -\beta}]\ . 
\]
In fact, $\kappa_{\alpha}$ corresponds to a  joint cumulant between 
$$(x,\stackrel[\alpha_{1,1} \text{times}]{}{\underbrace{X_{1,1},\ldots,  X_{1,1}}}, \ldots, \stackrel[\alpha_{n,n} \text{times}]{}{\underbrace{X_{n,n},\ldots,  X_{n,n}}}),$$
see e.g.~\cite{novak2014three} for an introduction to mixed cumulants. 
Note that the matrix $\alpha\in \mathbb{N}^{n\times n}$ can be considered as the adjacency matrix of a multigraph on the set $[n]$ of nodes. With a slight abuse of notation, we sometimes refer to $\alpha$ as a multigraph.
We write $N(\alpha)=\{i:\sum_{j}\alpha_{i,j}>0\}\subset [n]$ the set of active nodes of $\alpha$.
The following lemma states that, for most $\alpha$, we have $\kappa_{\alpha}=0$.  

\begin{lem}\label{lem:null:cumulant}
We have $\kappa_{\alpha}=0$ if either $\alpha$ contains self edges, that is $\max_{i}\alpha_{i,i}>0$, or if  $\alpha$ contains at least one non-trivial connected component that contains  neither the node $1$ nor the node $2$. 
\end{lem}
The proof of this lemma relies on the following fundamental property of the mixed cumulant $\kappa(Z_1,\ldots, Z_k)$ between random variables $Z_1,\ldots, Z_k$. If there exist a partition $(A,B)$ of $[k]$ such $(Z_i,i\in A)$ is independent of $(Z_i,i\in B)$, then $\kappa(Z_1,\ldots, Z_k)=0$.  See the proof of \Cref{lem:null:cumulant} for more details. 

\begin{proof}[Proof of Lemma~\ref{lem:null:cumulant}]
First, the random variables $X_{i,j}$ are constant and equal to $\lambda$. As a constant random variable is independent of any other random vector, we deduce that $\kappa_{\alpha}=0$ if $\max_{i}\alpha_{i,i}>0$. 

Let us turn to the second result and let us first assume that $\alpha$ has a non-trivial connected component that neither contains $1$ or $2$. 
We Write $N_1\subset N(\alpha)$ the set of nodes corresponding to these nodes and $N_2= (N(\alpha)\setminus N_1)\cup {\{1,2\}}$. For $s=1,2$, we observe $(X_{i,j})$, $(i,j)$ in $N_s$ is measurable with respect to the latent positions $(V_i)$, $i\in N_s$. Since the latent positions are independent, this enforces that $(X_{i,j})_{i,j\in N_1}$ is independent of $(X_{i,j})_{i,j\in N_2}$. Hence, by the fundamental property of cumulants, we have $\kappa_{\alpha}=0$.

\end{proof}

Hence, in~\eqref{eq:upper:low_correlation}, we only have consider two collections of multisets $\alpha$ without self edges.:
\begin{itemize}
    \item[(A)] Those such that $\{1,2\}\subset N(\alpha)$ and each connected component of $\alpha$ contains either $1$ or $2$.  
    \item[(B)] Those such that $ |N(\alpha) \cap \{1,2\} |= 1$ and $\alpha$ is connected. \end{itemize}

Then, we deduce from~\eqref{eq:upper:low_correlation} that 
\[
\mathrm{Corr}_{\leq D}^2 \leq \lambda^2 p_0^2 + \mathrm{Corr}_{\leq D,A}^2+ \mathrm{Corr}_{\leq D,B}^2\ ,  
\]
where 
\begin{equation*}
\mathrm{Corr}_{\leq D,A}^2 = \sum_{\alpha \in \mathbb{N}^N \, 1\leq |\alpha|\leq D,\  \alpha\text{ of type A}}\frac{\kappa_{\alpha}^2}{\alpha!} \ ; \quad 
\mathrm{Corr}_{\leq D,B}^2 = \sum_{\alpha \in \mathbb{N}^N \, 1\leq |\alpha|\leq D,\  \alpha\text{ of type B}}\frac{\kappa_{\alpha}^2}{\alpha!}\  . 
\end{equation*}

For the graphs of type $(A)$ a rough bound of the cumulant will be sufficient for our purpose. For the graphs of type (B), we need to prove that the corresponding cumulants are small enough. 

\begin{lem}\label{lem: correlation} Define $\rho = 2k/n$. Define $r_0 = 2\lambda^2(D+1)^4$. If $r_0 <1$ and $n\rho^2<1/2$, we have 
    \begin{eqnarray*}
        \mathrm{Corr}_{\leq D,A}^2&\leq& \lambda^2\rho^2\left(1 + \frac{4r_0}{1-r_0}\right) ; \quad 
        \mathrm{Corr}_{\leq D,B}^2\leq  \lambda^2 n\rho^4 \frac{r_0}{1-r_0} .
    \end{eqnarray*}    
\end{lem}

We split the proof of~\Cref{lem: correlation} and the control of $\mathrm{Corr}_{\leq D,A}^2$ and the one of $\mathrm{Corr}_{\leq D,B}^2$ into the next two subsections. From this lemma, we conclude that 
\[
\mathrm{MMSE}'_{\leq D}\geq  \lambda^2 p_0 - \lambda^2 \rho^2 \left(2 + \frac{5r_0}{1-r_0}\right).    
\]


\subsection{Control of the A-term} 
First we bound the number of multigraphs of type (A). 
\begin{lem}\label{lem:multigraph}
    For integers $d\geq 1$ and $2\leq h \leq d+1$, the number multigraphs $\alpha$ on $[n]$ such that (i) $|\alpha|=d$, (ii) $\{1,2\}\subset N(\alpha)$, (iii) $|N(\alpha)|=h$, (iv) each connected component contains either $1$ or $2$,  is at most 
    $2^{d}n^{h-2}h^{2d-h+2}$ 
\end{lem}

\begin{proof}[Proof of Lemma~\ref{lem:multigraph}]
    We can choose freely at most $n^{h-2}$ nodes. Since neither of these $h-2$ nodes are isolated and since they belong to a connected component of $1$ or $2$, there are at most $h^{h-2}$ choices of edges for connecting them to another one. Finally, we have $2h^2$ possibilities for each of the remaining $d-h+2$ edges. 
\end{proof}

\begin{lem}\label{lem:cumulant}
We have $\kappa_{0}\leq  \rho$. For  any $\alpha$ such that $|\alpha|\geq 1$ and $\alpha$ is either of type $A$ or of type $B$, we have 
\[
|\kappa_{\alpha}|\leq \lambda^{|\alpha|+1}(|\alpha|+1)^{|\alpha|} \rho^{|N(\alpha)|-1}
\] 
Besides, the only $\alpha$ such that $|\alpha|=1$ and $\alpha$ is of type $A$ is $\alpha=(1,2)$ and satisfies, $|\kappa_{\alpha}|\leq \lambda^2 \rho$. 
\end{lem}
Delaying the proof of \lemref{cumulant} to the end of this subsection and putting these two lemmata together, we conclude that 
\begin{eqnarray*}
    \mathrm{Corr}_{\leq D,A}^2&\leq& \lambda^2 \rho^2+  \sum_{d=2}^D \sum_{h=2}^{d+1} 2^{d}\lambda^{2(d+1)}(d+1)^{2d}n^{h-2}h^{2d-h+2}\rho^{2h-2} \\
    &\stackrel{(i)}{\leq} & \lambda^2 \rho^2\left(1+ 4 \sum_{d=2}^D\lambda^{2d}2^d(d+1)^{4d} \right) \\
    &\leq &  \lambda^2 \rho^2\left(1+ 4 \sum_{d=2}^D [2\lambda^{2}(d+1)^4]^{d} \right) \\ 
       &\leq & \lambda^2\rho^2\left(1 + \frac{2r_0}{1-r_0}\right)  , 
\end{eqnarray*}
where we used in (i) that $n\rho^2 \leq 1/2$ and $\lambda<1$ and, in the last line, that $r_0= 2\lambda^2(D+1)^4 < 1$.

\begin{proof}[Proof of Lemma~\ref{lem:cumulant}]
Denote $\kappa'_{\alpha}$ the cumulant $\kappa_{\alpha}$ in the specific case where $\lambda=1$. By multilinearity, we have $|\kappa_{\alpha}|= \lambda^{|\alpha|+1}\kappa'_{\alpha}$, so that we only have to focus on the case where $\lambda=1$, which we assume henceforth. We first consider the only $\alpha$ of type $A$ such that $|\alpha|=1$. One can readily check that is is equal to $\alpha=(1,2)$. For such $\alpha$, we have $\kappa_{\alpha}= \E[x^2]- \E[x]^2 = p_0-p_0^2 \leq \rho$.

Next, we prove the general bound \Cref{lem:cumulant} by induction on $|\alpha|$.
The bound is obviously true for $|\alpha|=0$ since no such multigraph is of type $(A)$ or type $(B)$. Denote $\#CC(\gamma)$ the number of connected components of $\gamma$. We claim that 
\begin{equation}\label{eq:bound_moment1}
\mathbb{E}[X^{\gamma}]\leq \rho^{|N(\gamma)|-\#CC(\gamma)}\ , \quad \quad \mathbb{E}[xX^{\gamma}]\leq \rho^{|N(\gamma'')|-\#CC(\gamma'')}\ , 
\end{equation}
where here, $\gamma''$ is the graph $\gamma$ where we have added the edge $(1,2)$. Let us prove~\eqref{eq:bound_moment1}. Note that the second bound is consequence of the first one.  Also, since the entries of $X$ are  either $0$ or $1$, we deduce that $\mathbb{E}[X^{\gamma}]\leq \mathbb{E}[X^{\gamma'}]$
where $\gamma'$ is a covering forest of $\gamma$. Since $\#CC(\gamma)=\#CC(\gamma')$, we only have to prove~\eqref{eq:bound_moment1} for forests $\gamma$. Again, by independence of the latent positions, we can restrict ourselves to the case where $\gamma$ is a tree. Then, we enumerate the nodes $N(\gamma)$ in such a way that, except for the first node, all arriving nodes are connected to a preceding one. Given the latent position of a node $i$, the probability that $X_{i,j}=1$ is smaller than $1$ almost surely. This implies~\eqref{eq:bound_moment1}.

Recall the recursive formula of cumulants.
\begin{equation}\label{eq:definition_kappa_2}
    \kappa_{\alpha}= \mathbb{E}[xX^{\alpha}]-\sum_{0\leq \beta \lneq \alpha }\kappa_{\beta}\binom{\alpha}{\beta}\mathbb{E}[X^{\alpha -\beta}]  . 
\end{equation}

Now, consider any $\alpha$ either of type $(A)$ or of type $(B)$. In the above equation~\eqref{eq:definition_kappa_2}, we can reduce our attention on the terms $\beta$ that are of type $(A)$ or $(B)$ otherwise the corresponding cumulant $\kappa_{\beta}$ is equal to zero. By assumption on $\alpha$ and by~\eqref{eq:bound_moment1}, we have $\mathbb{E}[xX^{\alpha}]\leq \rho^{|N(\alpha)|-1}$. We have also $\kappa_{0}\mathbb{E}[X^{\alpha}]\leq \rho^{|N(\alpha)|-1}$. For  any $\alpha$ and $\beta\neq \alpha$, we claim that --- the proof is provided below --- 
\begin{equation}\label{eq:upper:graph}
    |N(\beta)|+|N(\alpha-\beta)|-\#CC(\alpha -\beta)-\#CC(\beta) \geq |N(\alpha)|- \#CC(\alpha) \ . 
\end{equation}
Since $\#CC(\alpha)\leq 2$, this implies that 
\[
\rho^{|N(\beta)|+|N(\alpha-\beta)|-\#CC(\alpha -\beta)-1} \leq \rho^{|N(\alpha)|-1}\ .
\]
By induction hypothesis, we derive that 
\[
    \kappa_{\alpha}\leq \rho^{|N(\alpha)|-1}\left[1+ \sum_{0\leq \beta \lneq \alpha } (|\beta| +1)^{|\beta|} \right]\leq \rho^{|N(\alpha)|-1}[1+ (2^{|\alpha|}-1)|\alpha|^{|\alpha|-1}]\leq  \rho^{|N(\alpha)|-1} [|\alpha|+1]^{|\alpha|}\ . 
\]
This concludes the proof. 
\end{proof}

\begin{proof}[Proof of~\eqref{eq:upper:graph}]
Note that we do not change the right-hand side of~\eqref{eq:upper:graph} if we replace $\alpha$ by a covering forest of $\alpha$. In contrast, this cannot increase the left hand side term since for fixed $\gamma \leq \gamma'$, $|N(\gamma)| - |N(\gamma')| \geq \#CC(\gamma) - \#CC(\gamma')$. Without loss of generality we can therefore assume that $\alpha$ corresponds to a simple forest so that $|N(\alpha)|- \#CC(\alpha)$ stands for its number of edges $|\alpha|$. A a consequence, both $\beta$ and $\alpha-\beta$  are also forests so that $ |N(\beta)|+|N(\alpha-\beta)|-\#CC(\alpha -\beta)-\#CC(\beta)$ is also the number of edges of $\alpha$ --- since $|N(\beta)|-\#CC(\beta)$ is $|\beta|$ and $|N(\alpha - \beta)|-\#CC(\alpha - \beta)$ is $|\alpha - \beta|$, and $\beta \leq \alpha$. 
\end{proof}

\subsection{Control of the B-term} 
The number of multigraphs of type (B) as well as the corresponding cumulants  are bounded in the following lemmata
\begin{lem}\label{lem:multigraph2}
    For integers $d\geq 1$ and $2\leq h \leq d+1$, the number of connected multigraphs $\alpha$ on $[n]$ such that (i) $|\alpha|=d$, (ii) $|\{1,2\}\cap  N(\alpha)|=1$, (iii) $|N(\alpha)|=h$ is at most 
    $2^{d+1}n^{h-1}h^{2d-h+1}$ 
\end{lem}

\begin{proof}[Proof of Lemma~\ref{lem:multigraph2}]
Here, we can choose freely $h-1$ nodes and we need to choose one node among either $\{1\}$ or $\{2\}$. Since the graph is connected, there are $(2h)^{h-1}$ possibilities to add edges. Then, we have $2^{d-h+1}h^{2(d-h+1)}$ possibilities for adding the remaing edges. 
\end{proof}

\begin{lem}\label{lem:cumulant2}
 For  any $\alpha$ such that $\alpha$ is of type $B$, we have 
    \[
    |\kappa_{\alpha}|\leq \lambda^{|\alpha|+1}(|\alpha|+1)^{|\alpha|} \rho^{|N(\alpha)|}  
    \] 
\end{lem}
Putting these two lemmata together we conclude that 
\begin{eqnarray*}
    \mathrm{Corr}_{\leq D,B}^2&\leq& \sum_{d=1}^D \sum_{h=2}^{d+1} 2^{d+1}\lambda^{2(d+1)} n^{h-1}h^{2d-h+1} (d+1)^{2d}\rho^{2h} \\
    &\stackrel{(i)}{\leq} & 4n\rho^4\sum_{d=1}^D  2^{d}\lambda^{2(d+1)}  (d+1)^{4d}   \\
    &\leq &  4\lambda^2n \rho^4\sum_{d=1}^D [2\lambda^2 (D+1)^4]^d  \\ 
       &\leq & \lambda^2 n\rho^4 \frac{r_0}{1-r_0}  . 
\end{eqnarray*}
where we used in (i) that $n\rho^2 \leq 1/2$ and in the last line that $r_0= 2\lambda^2(D+1)^4 < 1$.


\begin{proof}[Proof of Lemma~\ref{lem:cumulant2}]
    As for the proof of Lemma~\ref{lem:cumulant2}, we only have to consider the case $\lambda= 1$.
Also, as in that proof,  we argue by induction on $\alpha$. If $|\alpha|=1$, we have $|\kappa_{\alpha}|\leq \mathbb{E}[xX^{\alpha}]+ \mathbb{E}[x]\mathbb{E}[X^{\alpha}] \leq 2\rho^2$. Now assume that $|\alpha|\geq 2$. 
By the formula \eqref{eq:definition_kappa_2}, we have 
\[
    |\kappa_{\alpha}|\leq |\mathbb{E}[xX^{\alpha}]|+ \sum_{0\leq \beta \lneq \alpha }|\kappa_{\beta}\mathbb{E}[X^{\alpha -\beta}]|\ ,    
\]
where the sum over $\beta$ only runs over $\beta$ that are of type (B), since $\beta$ cannot be of type (A) and otherwise $\kappa_{\beta}=0$. Relying on~\eqref{eq:upper:graph} and on the induction hypothesis we derive that 
\begin{eqnarray*}
    |\kappa_{\alpha}|&\leq& \rho^{|N(\alpha)|} +  \sum_{0\leq \beta \lneq \alpha }(|\beta|+1)^{|\beta|} \rho^{|V(\beta)|+|V(\alpha- \beta)|-\#CC(\alpha-\beta) }\\
&\leq & \rho^{|N(\alpha)|}\left[1 +  \sum_{0\leq \beta \lneq \alpha }(|\beta|+1)^{|\beta|}\right]\leq \rho^{|N(\alpha)|}(|\alpha|+1)^{|\alpha|}\  . 
\end{eqnarray*}
 \end{proof}

\section{Proofs of Subsection~\ref{sec:noncomp}}\label{proof:noncomp}

\subsection{Proof of Theorem~\ref{thm:noncompl2}}

For any $\Pi \in \cS_n$ and $\tilde \theta\in [0,A]^n~\mathrm{non-increasing}$ we have that
$$\|\Pi\cdot Y - T(\tilde \theta)\|_F^2 = \|Y - \Pi^\top\cdot T(\tilde \theta)\|_F^2.$$
So that
$$\|\Pi \cdot  Y - T(\tilde \theta)\|_F^2 =\|X - \Pi^\top\cdot  T(\tilde \theta)\|_F^2 + 2 \langle X - \Pi^\top \cdot T(\tilde \theta), E\rangle + \|E\|_2^2.$$

Since the entries of $E$ are i.i.d.~and distributed as $\cN(0,1)$ we have that
$$\langle X - \Pi^\top \cdot T(\tilde \theta), E\rangle \sim \cN(0, \|X - \Pi^\top \cdot T(\tilde \theta)\|^2_F).$$
So that for any $1/2>\delta>0$, by an union bound and properties of Gaussian tails, with probability larger than $1-\delta$: $\forall \Pi \in \cS_n, \forall \tilde \theta \in \cA^n$,$$\frac{\Big|\langle X - \Pi^\top \cdot T(\tilde \theta), E\rangle\Big|}{\|X - \Pi^\top \cdot  T(\tilde \theta)\|_F} \leq 2\sqrt{\log\Big(\frac{n! |\cA|}{\delta}\Big)} \leq 2\sqrt{n\log\big(n^3A\big)} + 2\sqrt{\log\Big(\frac{1}{\delta}\Big)},$$
since $|\cA^n|\leq \big(A n^2\big)^n$. So that on an event $\xi$ of probability larger than $1-\delta$: $\forall \Pi \in \cS_n, \forall \tilde \theta \in \cA^n$,
$$\Big|\|\Pi \cdot Y - T(\tilde \theta)\|_F^2 - \|\Pi \cdot X - T(\tilde \theta)\|_F^2 - \|E\|_F^2\Big| \leq  d\|X - \Pi^\top \cdot T(\tilde \theta)\|_F,$$
where $d = \Big[2\sqrt{n\log\big(n^3A\big)} + 2\sqrt{\log\Big(\frac{1}{\delta}\Big)}\Big]$.

Note that if $\inf_{\tilde \theta\in \cA^n~\mathrm{non-increasing}} \|X -  \Pi^\top \cdot T(\tilde \theta)\|_F\geq 2d$, then on $\xi$ we have $\inf_{\tilde \theta\in \cA^n~\mathrm{non-increasing}}\|\Pi \cdot Y - T(\tilde \theta)\|_F^2 - \|E\|_F^2 \geq d^2.$ 

Note also that by definition of $\cA^n$, there exists $\theta'\in \cA^n$ such that $\|\theta' - \theta\|_2 \leq u\sqrt{n}$, so that $\|T(\theta) - T(\theta')\|_F \leq nu = 1/n$. This implies in particular that
$$\|X - (\Pi^*)^\top\cdot  T(\theta')\|_F \leq 1/n,$$
as $\|X - (\Pi^*)^\top \cdot T(\theta)\|_F = 0$. So that on $\xi$
$$\Big|\inf_{\tilde \theta\in \cA^n~\mathrm{non-increasing}}\|\Pi^* \cdot Y - T(\tilde \theta)\|_F^2 - \|E\|_F^2\Big| \leq d/n \leq d^2/2.$$

We conclude from this that on $\xi$, an event of probability larger than $1-\delta$
$$\inf_{\tilde \theta\in \cA^n~\mathrm{non-increasing}} \|X -  (\hat \Pi^{(2)})^\top \cdot  T(\tilde \theta)\|_F \leq 2d,$$
so that
$$ \inf_{\tilde \theta\in [0,A]^n~\mathrm{non-increasing}} \|\hat \Pi^{(2)} \cdot  X -   T(\tilde \theta)\|_F = \inf_{\tilde \theta\in [0,A]^n~\mathrm{non-increasing}} \|X -  (\hat \Pi^{(2)})^\top\cdot  T(\tilde \theta)\|_F \leq 2d.$$
This concludes the proof as over the grid the loss is always bounded by $A^2n^2$.

\subsection{Proof of Theorem~\ref{thm:noncompl2la}}

Set $m = n^8$. Write $\pi_m(.)$ the projection of $v \in \mathbb R^+$ on the largest element of $\mathbb N/m$ smaller than $v$. Note that for $v, v'\in[0,1]$, we have
$$v - v'- 1/m \leq \pi_m(v) - \pi_m(v')\leq v - v'+ 1/m.$$

Consider $G_m = [-1,1] \cap (\mathbb Z/m)$. The function $\varphi$ is unimodal and takes value in $[0,A]$. So that the set
$$\mathcal U_m = \{u \in G_m: |\phi(u) - \phi(u+1/m)|> A/\sqrt{m}\},$$
is such that $|\mathcal U_m| \leq 2\sqrt{m}$.

Set $\mathcal X_m = [-1,1] \setminus \{[u-2/m,u+2/m), u \in \mathcal U_m\}$. We therefore have that for any $x\in[0,1]$
$$\lambda(\mathcal X_m \cap [-x,1-x]) \geq 1 - 32/\sqrt{m}.$$
So that for $X\sim \mathrm{Unif}([0,1])$. 
$$\mathbb P(X-x \in \mathcal X_m) \geq 1 - 32/\sqrt{m}.$$
So by a union bound we have that 
$$\mathbb P(\forall i,j, V_i - V_j \in \mathcal X_m) \geq 1 - n^2\times 32/\sqrt{m} = 1 - 32/n^2.$$
Let us write $\xi$ for the event where the above bound holds.

Note also that for any $v,v' \in [0,1]$ such that $v-v' \in \mathcal X_m$,  we have since $|\pi_m(v)- \pi_m(v') - v+v'| \leq 2/m$
$$|\varphi(\pi_m(v)- \pi_m(v')) - \phi(v-v')|\leq 2A/\sqrt{m} = 2A/n^4.$$
So on $\xi$, we have, writing $\tilde X_{i,j} = \pi_m(\varphi(\pi_m(V_i) - \pi_m(V_j) ))$ that
$$\|X - \tilde X\|_F^2 \leq 3A/n^2.$$
So that on $\xi$, we have that there exist $\tilde V \in \mathcal V^n$ non-decreasing and $\phi \in \mathcal A^n$ non-increasing such that
$$\|\Pi^*.X - R(\tilde V, \phi)\|_F^2 \leq 3A/n^2.$$
Using this and a similar proof as in the proof of Theorem~\ref{thm:noncompl2}, we conclude the proof.

\subsection{Proof of Theorem~\ref{thm:LBl2}}

Write
$$\cX = \{\Pi \cdot T(\theta), \Pi \in \cS_n\}.$$
Write $S_n$ for the set of functions from $\{1,\ldots,n\}$, namely the set of all bijections from $\{1,\ldots,n\}$ to $\{1,\ldots,n\}$. Write for any $\pi \in S_n$ and any $i \in \{1,\ldots,n\}$ 
$$f_\pi(i) =(\pi(i-1), \pi(i), \pi(i+1))$$
Write $d(\pi,\pi') = \sum_{i=1}^n \mathbf 1\{\exists j: f_\pi(i) = f_{\pi'}(j)~\mathrm{or}~f_{\tau\pi}(i) = f_{\pi'}(j)\}$.

\begin{lem}[Packing set of permutations in $d$-distance]
There exist two universal constants $c, c'>0$ and a set $\bar S_n \subset S_n$ such that:
\begin{itemize}
    \item $|\bar S_n| \geq (1+c)^n$.
    \item for any $\pi,\pi' \in \bar S_n$, $d(\pi,\pi') \geq c'n$.
\end{itemize}
\end{lem}
\begin{proof}
    This is a direct consequence of Gilbert-Varshamov bound, see e.g.~\cite[Lemma 15.3]{wu2017lecture}, by taking a grid of $i$ separated by a distance of $2$ and by lower bounding the number of partitions in triplets associated to these $i$ that differ in at least a fraction of the groups.
\end{proof}

The following corollary follows immediately.
\begin{cor}\label{cor:pack}
    There exists a universal constant $c>0$ and a set $\bar \cX \subset \cX$ such that:
\begin{itemize}
    \item $|\bar \cX| \geq (1+c)^n$.
    \item for any $X,X' \in \bar \cX$, $\|X-X'\|_F^2 \geq c' u^2n$.
\end{itemize}
\end{cor}
\begin{proof}
For any $X,X' \in \cX$ corresponding to two permutations $\pi,\pi' \in S_n$, we have
$$d(\pi,\pi') \leq u^{-2}\|X - X'\|_F^2 .$$
The corollary follows directly, with $\bar \cX$ taken as the set of permuted matrices defined through $\bar S_n$, namely the set of matrices $X$ such that there exists $\pi \in \bar S_n$ such that $X$ permuted through $\pi$ is equal to $T(\theta)$.
\end{proof}

Write for any $X\in \cX$, $\bbP_X$ for the distribution of $Y = X+E$, where $E$ is such that the $E_{i,j}$ are i.i.d.~$\mathcal N(0,1)$. We have for any $X\in \cX$ that
\begin{equation}\label{eq:KL}
    \mathrm{KL}(\bbP_X, \bbP_{T(\theta)}) = \frac{\|X - T(\theta)\|_F^2}{2} \leq 3u^2n/2,
\end{equation}
where $\mathrm{KL}(.,.)$ is the Kullback-Leibler divergence.

We now remind the well-known Fano's lemma, see e.g.~\cite[Corollary 2.6]{tsybakov2003introduction}.
\begin{lem}[Fano's lemma]
    Consider, for $m\geq 2$, $m+1$ probability distributions $P_0, P_1,\ldots, P_m$, and an estimator $\psi$ based on a sample from $P \in \{P_0, P_1,\ldots, P_m\}$ that takes value in $\{1,\ldots, m\}$. Then
    $$\max_{j\leq m} P_j(\psi \neq j)  \geq 1 - \frac{ \log(2) + 2\max_{1\leq j\leq m} \mathrm{KL}(P_j,P_0)}{\log(m)}.$$
\end{lem}
We apply Fano's lemma $\bar \cX$ from Corollary~\ref{cor:pack}. We have for any measurable function $\psi$ of the data
$$\max_{X\in \bar \cX} \bbP_X(\psi \neq X) \geq 1 - \frac{ \log(2) + 2\max_{X\in \bar \cX} \mathrm{KL}(\bbP_X,\bbP_{T(\theta)})}{\log(|\bar \cX|)} \geq 1 - \frac{ \log(2) + 2\max_{X\in \bar \cX} \mathrm{KL}(\bbP_X,\bbP_{T(\theta)})}{n\log(1+c)}.$$
As for any $X'X' \in \bar \cX$ such that $X\neq X'$ we have $\|X - X'\|_F^2 \geq nu^2/4$ by Corollary~\ref{cor:pack}. This implies for any estimator $\hat \Pi \in \cS_n$, by triangular inequality and Markov inequality
$$\frac{1}{(c')^2nu^2}\max_{X\in \bar \cX} \bbE_X \|\hat \Pi\cdot X - T(\theta)\|_F^2 \geq 1 - \frac{ \log(2) + 2\max_{X\in \bar \cX} \mathrm{KL}(\bbP_X,\bbP_{T(\theta)})}{n\log(1+c)}.$$
By Equation~\eqref{eq:KL}
$$\frac{1}{(c')^2nu^2}\max_{X\in \bar \cX} \bbE_X \|\hat \Pi\cdot X - T(\theta)\|_F^2 \geq 1 - \frac{ \log(2) + 3u^2n}{n\log(1+c)}.$$
This concludes the proof as for $X \in \cX$, we straightforwardly have for any $\hat \Pi$:
$$\inf_{R \in \mathcal{R}_n} \|\hat \Pi\cdot X -   R\|_F \geq \frac{1}{2}\|\hat \Pi\cdot X -   T(\theta)\|_F.$$

\subsection{Proof of Theorem~\ref{thm:LBl2l}}


Consider a Poisson random variable $\tilde n \sim \mathcal P(n)$. 
Consider in what follows $\tilde V_1, \ldots, \tilde V_{\tilde n} \sim_{\mathrm{i.i.d.}} \mathrm{Unif}([0,1])$, and we consider the matrix $\tilde X$ of size $\tilde n \times \tilde n$, associated with $\tilde V_1, \ldots, \tilde V_{\tilde n}$.


Define $\mathcal G = [0,1] \cap (6a\mathbb N)$. Consider $x \in \mathcal G$. Define
$$Z_x = \mathbf 1\{\exists i,j \leq \tilde n: |\tilde V_i - x|\lor |\tilde V_j - x| \leq a/2,~~\text{and}~~\forall k\not \in \{i,j\}, |\tilde V_k - x| \geq 3a\}.$$
Since $\tilde n$ is a Poisson random variable, we have that the $(Z_x)_{x \in \mathcal G}$ are i.i.d.~and are Bernoulli random variables with parameter
$$p = 2^{-1}(an/2)^2 \exp(-an) \times \exp(- 6an) = \exp(-7)/8,$$
since $a = 1/n$.

Consider now $\mathcal U = \{x\in \mathcal G: Z_x = 1\},$ and write for any $x\in \mathcal U$ resp.~$\{i_x,j_x\}$ for the indexes for the two points $\tilde V$ that are in $[x-a/2, x+a/2]$. Note that for $x\in \mathcal U$, the rows $\tilde X_{i_x}, \tilde X_{j_x}$ are such that they take value $0$ everywhere but in $i_x,j_x$, where they take value $u$.

Assume then that an oracle gives us the set of indexes
$$\{i_x, j_x, x\in \mathcal U\},$$
and that our aim now is to recover the pairs $\{\{i_x, j_x\}, x \in \mathcal U\}$. This aim is simpler as seriation, in the sense that
$$\inf_{\hat \Pi} \sup_\Pi \mathbb E_\Pi [\ell_2(\hat \Pi,\Pi)] \geq \inf_{\{\{\hat i_x, \hat j_x\}, x \in \mathcal U\}} \sup_\Pi \mathbb E_\Pi \left[ u^2 d(\{\{\hat i_x, \hat j_x\}, x \in \mathcal U\}, \{\{ i_x,  j_x\}, x \in \mathcal U\})\right],$$
where $\{\{\hat i_x, \hat j_x\}, x \in \mathcal U\}$ forms a partition of $\{i_x, j_x, x\in \mathcal U\}$, and where $d(\{\{\hat i_x, \hat j_x\}, x \in \mathcal U\}, \{\{ i_x,  j_x\}, x \in \mathcal U\})$ is the number of pairs which differ between the two partitions.

Set $d = |\mathcal U|$. So from now on, we consider the simpler problem, which is akin to seriate a dimension $2d \times 2d$ matrix $X'$ such that
$$X_{i_k, j_k}' = X_{j_k, i_k}' = X_{i_k,i_k}' = X_{i_k,i_k}' =  u,~~~\forall k \leq d,$$
and otherwise it takes value $0$, and where $\{\{i_k,j_k\}, k \leq d\}$ are pairs belonging to $\{1,\ldots,2d\}^2$, which together form a partition of $\{1,\ldots,2d\}$. The seriation problem is here equivalent to recover the pairs, and the error we obtain is higher than
$$u^2 d(\{\{i_k,j_k\}, k \leq d\}, \{\{\hat i_k,\hat j_k\}, k \leq d\}),$$
where $\{\{\hat i_k,\hat j_k\}, k \leq d\}$ is the estimated partition in pairs.

We the proceed similar to the proof of Theorem~\ref{thm:LBl2}:
\begin{itemize}
    \item We first prove that there exists a packing set of the set of all partitions of $\{1, \ldots, d\}$ in pairs, with distance $c'd$, and which cardinality larger than $(1+c)^d$ for some universal constant $c>0$.
    \item We prove that the KL between any two matrices $X'$ as described above is upper bounded by $4u^2 d$.
    \item We apply Fano's lemma and obtain for any vector $v\in [0,1]^n$
    $$\inf_{\hat \Pi} \sup_{\Pi} \mathbb E [\ell_2(\hat \Pi, \Pi)|V=v] \geq \Big[1 - \frac{\log(2)+8u^2d}{d(\log(1+c))}\Big] \frac{du^2}{16}.$$
\end{itemize}
Since by definition of $d = |\mathcal U|$, and since by definition $|\mathcal U| = \sum_{x\in \mathcal G} Z_x$, we have that with probability larger than $1-1/n^4$
$$d\geq |\mathcal G| p - 4\sqrt{|\mathcal G| \log(n)} \succeq n.$$
So that the result follows, also by concentration of the Poisson random variable $\tilde n$ with exponential probability in $[n/2,2n]$.

\section{Proof of Proposition~\ref{prp:LBlinf}} \label{app:LBinf}
In what follows, consider $\theta =  (1,1,0,\ldots,0)$ and
$$\cX = \{\Pi \cdot T(\theta), \Pi \in \cS_n\}.$$
    Now consider any two $X^{(1)},X^{(2)} \in \mathcal X$ such that $X^{(1)}\neq X^{(2)}$, and consider two different noisy datasets:
$$Y^{(1)} = X^{(1)} + E^{(1)},$$
and
$$Y^{(2)} = X^{(2)}+E^{(2)},$$
where $E^{(1)} = 0$ and $E^{(2)} = X^{(1)} - X^{(2)}.$ Note that $\|E^{(1)}\|_\infty = 0 \leq 1$ and $\|E^{(2)}\|_\infty \leq \|X^{(1)} - X^{(2)}\|_\infty \leq 1$, so that the adversarial noise satisfies our assumption. Note also that $Y^{(1)} = Y^{(2)}$, and $X^{(1)}\neq X^{(2)}$, so that for any estimator $\hat \Pi \in \cS_n$ of the permutation, we either have
$$\hat \Pi\cdot X^{(1)} \neq T(\theta),$$
or
$$\hat \Pi\cdot X^{(2)} \neq T(\theta).$$
By definition of $\cX$ and since $\hat \Pi \in \cS_n$, this implies that
$$\|\hat \Pi\cdot X^{(1)} -  T(\theta)\|_\infty \geq 1,$$
or
$$\|\hat \Pi\cdot X^{(2)} -  T(\theta)\|_\infty \geq 1,$$
which implies
$$\|\hat \Pi\cdot X^{(1)} -  T(\theta)\|_\infty \lor \|\hat \Pi\cdot X^{(2)} -  T(\theta)\|_\infty \geq 1.$$
This concludes the proof.

\section{A pseudo-code algorithm for \texttt{PINES}} \label{app:algo}

\begin{algorithm}[ht]
\caption{The \algoname algorighm}\label{algo:generic}
\KwData{a distance $\hat \dt$ on $[n] \times [n]$ and connectivity radii $\rho_1,\rho_2,\rho_3 > 0$}
\KwResult{a permutation $\Pi$ of $[n]$
}
$D \gets \emptyset$\;
\While{$D \neq [n]$ \tcc{Building the packing $\cP$ and the associated partition}}{
Take $i \in [n] \setminus D$\;
$\cP \gets \cP \cup \{i\}$\;
$P_i \gets \{j \in [n] \setminus D~|~\hat\dt(i,j) \leq \rho_1\}\setminus D$\;
$D \gets D \cup P_i$\;
}
\For{$i \in \cP$ \tcc{Build the neighborhood graphs and compute the connected components}}{
$V_i \gets \{j \in [n]~|~\hat \dt(i,j) > \rho_2\}$\;
$G_i \gets$ $\rho_3$-neighborhood graph on $V_i$\;
$\cC_i \gets$ the equivalence classes of $\cP \setminus \{i\}$\;
The equivalence relation is $j \sim k$ if and only if $j$ and $k$ are connected in $G_i$
}
$Q \gets \emptyset$\;
\tcc{Seriate the packing}
\eIf{there exists $i \in \cP$ such that $\Card \cC_i = 1$ \tcc{Find an extremal point}}
{$Q \gets \{i\}$\;
$\Pi_{\cP}(1) \gets i$\;
\While{$Q \neq \cP$}
{
\eIf{there exists $i \in \cP \setminus Q$ and $C \in \cC_i$ such that $Q = C$ \tcc{Find next point}}
{
Append $i$ to $Q$\;
$\Pi_{\cP}(\Card Q) \gets i$\;
}{\textbf{Raise Error}\;}
}
}{\textbf{Raise Error}\;}
$q \gets 1$ \tcc{A counter to know how many points we have ordered yet.}
\For{$r \in [\Card \cP]$}
{
$i \gets \Pi_\cP(r)$\;
\For{$k \in P_i$ \tcc{Arbitrarily order within $P_i$}}{
$\Pi(k) \gets q$\;
$q \gets q+1$\;
}
}
\textbf{return} $\Pi$.
\vspace{.5em}
\end{algorithm}



\end{document}